\documentclass[11pt]{article}
\usepackage[latin1]{inputenc}
\usepackage{float,epsfig}
\usepackage{amssymb}
\usepackage[algonl,algo2e,algoruled]{algorithm2e}
\usepackage{color}
\graphicspath{{./images/}}

\newcommand\eqdef{\stackrel{{\rm def}}{=}}

\renewcommand{\bar}[1]{\overline{#1}}
\renewcommand{\hat}[1]{\widehat{#1}}
\renewcommand{\tilde}[1]{\widetilde{#1}}

\def\e{\epsilon}

\def\rank{\mathrm{rank\,}}

\def\tr{\mathrm{trace\,}}

\newtheorem{lemma}{Lemma}
\newtheorem{theorem}{Theorem}
\newtheorem{proposition}{Proposition}

\newtheorem{assumption}{Assumption}

\def\E{\bold E}
\def\P{\bold P}

\newcommand{\eE}{\mathbb{E}}

\newcommand{\epr}{\hfill \( \Box \quad \)\\}
\newcommand{\R}{\mathbb{R}}

\newcommand{\cA}{\mathcal{A}}
\newcommand{\cB}{\mathcal{B}}

\newcommand{\cF}{\mathcal{F}}

\newcommand{\cI}{\mathcal{I}}

\newcommand{\cL}{\mathcal{L}}

\newcommand{\cN}{\mathcal{N}}
\newcommand{\cO}{\mathcal{O}}
\newcommand{\cP}{\mathcal{P}}

\newcommand{\cU}{\mathcal{U}}

\newcommand{\cW}{\mathcal{W}}
\newcommand{\cX}{\mathcal{X}}
\newcommand{\cY}{\mathcal{Y}}
\newcommand{\cZ}{\mathcal{Z}}

\newcommand{\argmin}{\mathrm{arg\,min}}

\newcommand{\diag}{\mathrm{diag}}

\newcommand{\be}{\begin{eqnarray}}
\newcommand{\ee}[1]{\label{eq:#1}\end{eqnarray}}
\newcommand{\nn}{\nonumber \\}
\newcommand{\ese}{\end{eqnarray*}}
\newcommand{\bse}{\begin{eqnarray*}}
\newcommand{\rf}[1]{~(\ref{eq:#1})}
\renewcommand{\(}{$\,}
\renewcommand{\)}{\,$}
\title{Sparse Non Gaussian Component Analysis by Semidefinite Programming}
\author{Elmar~Diederichs \thanks{The author is partially supported by the Laboratory of Structural
Methods of Data Analysis in Predictive Modeling, MIPT, through the
RF government grant, ag. 11.G34.31.0073. Furthermore the author is
grateful for funding to the Competitive Procedure of the Leibniz
Association within the ''Pact for Research and Innovation''
framework.}\\
Weierstrass Institute\\
Mohrenstr. 39, 10117 Berlin, Germany\\
{\tt diederic@wias-berlin.de}
\and Anatoli Juditsky
\\
LJK,
Universit\'e J. Fourier, BP 53 38041\\
Grenoble cedex 9, France \\
{\tt juditsky@imag.fr}
\and Arkadi Nemirovski\\
ISyE, Georgia Institute of Technology\\
 Atlanta, Georgia 30332, USA\\
 {\tt nemirovs@isye.gatech.edu}
\and Vladimir Spokoiny \\
Weierstrass Institute and Humboldt University,\\
Mohrenstr. 39, 10117 Berlin, Germany\\
{\tt spokoiny@wias-berlin.de}
}
\begin{document}

%\makeanontitle
\maketitle

\begin{abstract}
Sparse non-Gaussian component analysis (SNGCA) is an
unsupervised method of extracting a linear structure from a high
dimensional data based on estimating a low-dimensional non-Gaussian data
component. In this paper we discuss a new approach to direct estimation of
the projector on the target space  based on semidefinite programming which improves the method sensitivity to a broad variety of
deviations from normality.
\par
We also discuss the procedures which allows to recover the structure when its effective dimension is unknown.
\\[2ex]
{\noindent}{\em Mathematical Subject Classification:} 62G07, 62H25, 62H30,
90C48, 90C90
\end{abstract}

%-------------------------------------------------------------------------
\section{Introduction}
%-------------------------------------------------------------------------
Numerous statistical applications are confronted today with the
so-called curse of dimensionality (cf. \cite{HTF,W}). Using high-dimensional datasets implies an exponential increase of
computational effort for many statistical routines, while the data thin out in the local
neighborhood of any given point and classical statistical methods
become unreliable. When a random phenomenon is observed in the high
dimensional space \( \R^{d} \) the ''intrinsic dimension'' $m$ covering
degrees of freedom associated with same features may be much smaller than
$ d$. Then introducing \emph{structural assumptions} allows to reduce the
problem complexity without sacrificing any statistical information
\cite{M,RS}. In this study we consider the case where the phenomenon
of interest is (approximately) located in a linear subspace $\cI$.
When compared to other approaches which involve construction of nonlinear mappings
from the original data space onto the ''subspace of interest'', such that
isomaps \cite{TSL}, local-linear embedding \cite{HTF} or Laplacian
eigenmaps \cite{bn}, a linear mapping appears attractive due to its
simplicity --- it may be identified with a simple object, the  projector $\Pi^*$ from $\R^d$
onto $\cI$. To find the structure of interest a statistician may seek for
the non-Gaussian components of the data distribution, while its Gaussian
components, as usual in the statistical literature %$\cite{CT},
may be
treated as non-informative noise.

%Sparse non-gaussian component analysis (SNGCA), introduced in \cite{DJSS}, is
%a linear structure extraction method for
%semiparametric high dimensional data analysis which is based on estimating
%low-dimensional non-Gaussian components of the high-dimensional data distribution.
%SNGCA is structural adaptive in the sense that every new step essentially
%uses the result of previous iterations and does not assume any a priori
%knowledge about \(\rho\).
Several techniques of estimating the ``non-Gaussian subspace'' have been proposed
recently. In particular, NGCA (for {\em Non-Gaussian Component Analysis}) procedure,
introduced in \cite{BKSS}, and then developed into SNGCA (for {\em Sparse NGCA}) in
\cite{DJSS}, is based on the decomposition the problem of dimension reduction into two
tasks: the first one is to extract from the data a set $\{\widehat{\beta}_j\}$ of
candidate vectors \(\widehat{\beta}_j\) which are ''close'' to \( \cI \). The second is
to recover an estimation $\hat{\Pi}$ of the projector $\Pi^*$ on \( \cI \) from
$\{\widehat{\beta}_j\}$. In this paper we discuss a new method of SNGCA based on
Semidefinite Relaxation of a nonconvex minmax problem which allows for a direct recovery
of $\Pi^*$. When compared to previous implementations of the SNGCA in
\cite{BKSS,KawSugBlaMue07,DJSS}, the new approach ''shortcuts'' the intermediary stages
and makes the best use of available information for estimation of $\cI$. Furthermore, it
allows to treat in a transparent way the case of unknown dimension $m$ of the target
space $\cI$.

The paper is organized as follows: in Section~\ref{SFramework}
we present the setup of SNGCA and briefly review some existing techniques.
Then in Section~\ref{sec:ssdp} we introduce the new approach to recovery of
the non-Gaussian subspace and analyze its accuracy. Further we provide a
simulation study in Section~\ref{numerics}, where we compare the
performance of the proposed algorithm SNGCA to that of some known
projective methods of feature extraction.
%
%-------------------------------------------------------------------------
\section{Sparse Non-Gaussian Component Analysis}\label{SFramework}
%-------------------------------------------------------------------------

\subsection{The setup}
The \emph{Non-Gaussian Component Analysis} (NGCA)  approach is based on the assumption
that a high dimensional distribution tends to be normal in almost any randomly
selected direction.
This intuitive fact can be justified by the central limit theorem when the number of
directions tends to infinity.
It leads to the NGCA-assumption:
the data distribution is a superposition of a full dimensional Gaussian distribution and
a low dimensional non-Gaussian component.
In many practical problems like clustering or classification, the Gaussian component is
uninformative and it is treated as noise.
The approach suggests to identify the non-Gaussian component and to use it for the further
analysis.

The NGCA set-up can be formalized as follows; cf. \cite{BKSS}.
Let \( X_{1},...,X_{N} \) be i.i.d. from a
distribution \( \P \) in \( \R^{d} \) describing the random phenomenon
of interest. We suppose that \( \P \) possesses a density \( \rho \)
w.r.t. the Lebesgue measure on \( \R^{d} \), which can be decomposed as
follows:
\be
    \rho(x)=\phi_{\mu,\Sigma}(x)q(Tx).
\ee{modeldensity}
Here \( \phi_{\mu,\Sigma} \) denotes the density of
the multivariate normal distribution \( \cN(\mu,\Sigma) \) with parameters
\( \mu\in \R^{d} \) (expectation) and \( \Sigma\in \R^{d\times d} \)
positive definite (covariance matrix). The function \( q:\R^{m}\to \R \)
with \( m\leq d \) is positive and bounded. \( T\in \R^{m\times d}
\) is an unknown linear mapping. We refer to \( \cI={\rm
range}\; T \) as {\em target} or {\em non-Gaussian subspace}. Note that
though $T$ is not uniquely defined, \( \cI \) is well defined, same as the Euclidean projector \( \Pi^{*} \) on \( \cI \).% and can be
%considered as {\em the effective dimension reduction space} (EDR-space).
In what follows, unless it is explicitly specified otherwise, we assume that the {\em effective dimension} $m$ of
$\cI$ is known {\em a priori}. For the sake of simplicity we assume
that the expectation of $X$ vanishes: $\E[X]=0$.

The model \rf{modeldensity} allows for the following interpretation (cf. Section 2 of \cite{BKSS}): suppose that the observation $X\in \R^d$ can be decomposed into
$X = Z+\xi,$ where $Z$ is an ``informative low-dimensional signal'' such that $Z\in \cI$, $\cI$ being an $m$-dimensional subspace of $\R^d$, and $\xi$ is independent and Gaussian.
One can easily show (see, e.g., Lemma 1 of \cite{BKSS}) that in this case the
density of $X$ can be represented as \rf{modeldensity}.

\subsection{Basics of SNGCA estimation procedure} The estimation of $\cI$ relies upon
the following result, proved in \cite{BKSS}: suppose that the function $q$ is smooth, then for any smooth  function \(
\psi:\,\R^{d}\to \R \) the assumptions of \rf{modeldensity} and \( \E[X]=0
\) ensure that for
\be
    \beta(\psi)
    \eqdef
    \E \bigl[\nabla \psi(X)\bigr]
    =
    \int \nabla \psi(x) \, \rho(x)\, dx ,
\ee{originalbeta} there is  a vector \( \beta \in
\cI \) such that
\bse
    | \beta(\psi) - \beta |_{2}
    \leq
    \big| \Sigma^{-1} \E[ X \psi(X) ] \big|_{2}
\ese
where \( \nabla \psi \) denotes the gradient of \( \psi \) and $|\cdot|_p$ is the standard $\ell_p$-norm on $\R^d$. In particular,
if $\psi$ satisfies
\( \E[X \psi(X)]=0\), then \( \beta(\psi)\in \cI \). Consequently \be
    | (I - \Pi^{*}) \beta(\psi) |_{2}
    \le
    \Big| \Sigma^{-1} \int x \psi(x) \rho(x) \; dx \Big|_{2},
\ee{bounds} where \( I \) is the \( d \)-dimensional identity matrix and
\( \Pi^{*} \) is the Euclidean projector on \( \cI \).

The above result suggests the following two-stage estimation
procedure: first compute a set of estimates $\{\hat{\beta}_\ell\}$ of
elements \(\{\beta_j\}\) of \( \cI \), then recover an estimation of $\cI$
from $\{\hat{\beta}_\ell\}$. This heuristics has been first used to
estimate $\cI$ in \cite{BKSS}. To be more precise, the construction
implemented in \cite{BKSS} can be summarized as follows: let for a family
\( \{h_{\ell}\},\;\ell=1,...,L \) of smooth bounded (test) functions on \(
\R^{d} \) \be
 \gamma_{\ell} \eqdef \E[X h_{\ell}(X)],\;\; \eta_{\ell} \eqdef
\E[\nabla h_{\ell}(X)], \ee{glel} and let \be \hat{\gamma}_{\ell}\eqdef
N^{-1}\sum_{i=1}^N X_ih_{\ell}(X_i),\;\;\hat{\eta}_\ell \eqdef
N^{-1}\sum_{i=1}^N \nabla h_{\ell}(X_i) \ee{empir} be their ''empirical
counterparts''. The set of ''approximating vectors''
$\{\hat{\beta}_\ell\}$ used in \cite{BKSS} is as follows:
$\hat{\beta}_\ell=\hat{\eta}_\ell-\hat{\Sigma}^{-1}\hat{\gamma}_\ell,\;\ell=1,...,L$,
where $\hat{\Sigma}$ is an estimate of the covariance matrix $\Sigma$. The
projector estimation at the second stage is $\hat{\Pi}=\sum_{j=1}^m
e_je_j^T$, where $e_j,\;j=1,...,m$, are $m$ principal eigenvectors of the
matrix $\sum_{\ell=1}^L\hat{\beta}_\ell\hat{\beta}_\ell^T$. A numerical
study, provided in \cite{BKSS}, has shown that the above procedure can be
used successfully to recover $\cI$. On the other hand, such
implementation of the two-stage procedure possesses two important
drawbacks: it relies upon the estimation of the covariance matrix $\Sigma$
of the Gaussian component, which can be hard even for moderate dimensions
$d$. Poor estimation of $\Sigma$ then will result in
badly estimated vectors $\hat{\beta}_\ell$, and as a result, poorly
estimated $\cI$. Further, using the eigenvalue decomposition of the
matrix $\sum_{\ell=1}^L\hat{\beta}_\ell\hat{\beta}_\ell^T$ entails that
the variance of the estimation $\hat{\Pi}$ of the projector $\Pi^*$ on
$\cI$ is proportional to the number $L$ of test-functions.
As a result, the estimation procedure is restricted to utilizing only relatively small families
$\{h_\ell\}$, and is sensitive to the initial selection of ''informative'' test-functions.\\[2ex]

{\noindent}To circumvent the above limitations of the approach of
\cite{BKSS} a different estimation procedure has been proposed in
\cite{DJSS}. In that procedure the estimates $\hat{\beta}$ of vectors from
the target space  are obtained by the method, which was referred to as {\em
convex projection}. Let \( c \in \R^{L} \) and let
%\( \beta(c), \gamma(c) \in \R^{d} \) with
\begin{eqnarray*}
    \beta(c) = \sum_{l=1}^L c^{\ell} \eta_{\ell} \qquad
    \gamma(c) = \sum_{l=1}^L c^{\ell} \gamma_{\ell}.
\end{eqnarray*}
Observe that \( \beta(c) \in \cI \) conditioned that \(
\gamma(c)=0\). Indeed, if  \( \psi(x)=\sum_{\ell} c^{\ell} h_{\ell}(x) \),
then \( \sum_{\ell}c^{\ell}\E[ Xh_{\ell}(X)]=0 \), and by \rf{bounds},
\[
    \eta(c) = \sum_{\ell}c^{\ell}\E [\nabla h_{\ell}(X)] \in \cI.
\]
Therefore, the task of estimating $\beta\in\cI$ reduces to that of finding
a ''good'' corresponding coefficient vector. In \cite{DJSS} vectors
$\{\hat{c}_j\}$ are computed as follows: let
\begin{eqnarray*}
    \hat{\eta}(c)= \sum_{l=1}^L c^{\ell} \hat{\eta}_{\ell}\;\;\mbox{and}\;\;
    \hat{\gamma}(c)&=& \sum_{l=1}^L c^{\ell} \hat{\gamma}_{\ell}, \;\ell=1,...,L
\end{eqnarray*}
and let $\xi_j\in \R^{d},\;j=1,...,J$ constitute a set of \emph{probe}
unit vectors. Then it holds
\be
    \hat{c}_j
    =
    \argmin_{c}\left\{|\xi_j - \hat{\eta}(c) |_{2}\,\left|\;
    %\mbox{   subject to  }
    \hat{\gamma}(c)=0,\;| c |_{1} \le 1 \right.\right\},
\ee{chat}
and we set $\hat{\beta}_j=\hat{\beta}(\hat{c}_j)= \sum_{\ell} \hat{c}_j^\ell
\hat{\eta}_{\ell}$. Then $\cI$ is recovered by computing $m$ principal
axes of the minimal volume ellipsoid (Fritz-John ellipsoid) containing the
estimated points \( \{\pm\hat{\beta}_j\}_{j=1}^J \).

The recovery of $\hat{\cI}$ through the Fritz-John ellipsoid
(instead of eigenvalue decomposition of the matrix $\sum_\ell
\hat{\beta}_j\hat{\beta}_j^T$) allows to bound the estimation error of
$\cI$ by the maximal error of estimation $\hat{\beta}$ of elements of the
target space (cf. Theorem 3 of \cite{DJSS}), while the
$\ell_1$-constraint on the coefficients $\hat{c}_j$ allows to control
efficiently the maximal stochastic error of the estimations
$\hat{\beta}_j$ (cf. Theorem 1 of \cite{DJSS,Spexp2009}). On the other
hand, that construction heavily relies upon the choice of the probe vectors
$\xi_j$. Indeed, in order to  recover the projector on $\cI$, the collection of $\hat{\beta}_j$ should comprise at least $m$ vectors with non-vanishing
projection on the target space. To cope with this problem a multi-stage
procedure has been used in \cite{DJSS}: given a set $\{\xi_j\}_{k=0}$ of
probe vectors an estimation $\hat{\cI}_{k=0}$ is computed, which is used
to draw new probe vectors $\{\xi_j\}_{k=1}$ from the vicinity of
$\hat{\cI}_{k=0}$; these vectors are employed to compute a new estimation $\hat{\cI}_{k=1}$, and
so on. The iterative procedure improves significantly the accuracy of the
recovery of $\cI$. Nevertheless, the choice of ''informative'' probe
vectors at the first iteration $k=0$ remains a challenging task and
hitherto is a weak point of the procedure.

\section{Structural Analysis by Semidefinite Programming}
\label{sec:ssdp}
In the present paper we discuss a new choice of vectors $\beta$ which
solves the initialization problem of probe vectors for the SNGCA procedure
in quite a particular way. Namely, the estimation procedure we are to
present below does not require any probe vectors at all.

\subsection{Informative vectors in the target space}
Further developments are based on the following simple observation. Let
$\eta_\ell$ and $\gamma_\ell$ be defined as in \rf{glel}, and let
${U}=[\eta_1,...,\eta_L]\in \R^{d\times L}$, ${G}=[\gamma_1,...,\gamma_L
]\in \R^{d\times L}$. Using the observation in the previous section we
conclude that if $c\in \R^L$ satisfies
$Gc=\sum_{\ell=1}^Lc^\ell\gamma_\ell=0$ then
$Uc=\sum_{\ell=1}^Lc^\ell\eta_\ell$ belongs to $\cI$. In other words, if
$\Pi^*$ is the Euclidean projector on $\cI$, then \[(I-\Pi^*)Uc=0. \]
Suppose now that the set $\{h_\ell\}$ of test functions is
rich enough in the sense that vectors $Uc$ span $\cI$ when $c$ spans the
subspace  $Gc=0$. Recall that we assume the dimension $m$ of the
target space to be known. Then projector $\Pi^*$ on $\cI$ {\em is fully
identified} as the optimal solution to the problem \be
\Pi^*=\arg\min_\Pi\max_c\left\{|(I-\Pi)Uc|^2_2\left|\;
\begin{array}{c}
\mbox{$\Pi$ is a projector on an }\\\mbox{$m$-dimensional subspace of $\R^d$}\\
c\in \R^L, \;Gc=0
\end{array}\right.\right\}.
\ee{proj0} In practice vectors $\gamma_\ell$ and $\eta_\ell$ are not
available, but we can suppose that their ``empirical counterparts'' --
vectors $\hat{\gamma}_\ell,\;\hat{\eta}_\ell,\; \ell=1,...,L$ can be
computed, such that for a set $A$ of probability at least $1-\varepsilon$,
\be |\hat{\eta}_\ell-\eta_\ell|_2\le
\delta_N,\;\;|\hat{\gamma}_\ell-\gamma_\ell|_2\le \nu_N, \;\;\ell=1,...,L.
\ee{mdp} Indeed, it is well known (cf., e.g., Lemma 1 in \cite{DJSS} or
\cite{vdWaart96}) that if functions $h_{\ell}(x) = f(x,\omega_{\ell})$,
$\ell=1,...,L$, are used, where $f$ is continuously differentiable,
$\omega_\ell\in \R^d$ are vectors on the unit sphere and $f$ and $\nabla_x
f$ are bounded, then \rf{mdp} holds with
\be
\begin{array}{rcl}
\delta_N&=&C_1\max_{x\in \R^d,\,|\omega|_2=1}|\nabla_x f(x,\omega)|_2 N^{-1/2}\sqrt{\min\{d,\ln L\}+\ln\varepsilon^{-1}},\\
\nu_N&=&C_2 \max_{x\in \R^d,\,|\omega|_2=1} |xf(x,\omega)|_2
N^{-1/2}\sqrt{\min\{d,\ln L\}+\ln\varepsilon^{-1}},
\end{array}\ee{dnu} where $C_1,\,C_2$
are some absolute constants depending on the smoothness properties and the second moments
of the underlying density.

Then for any $c\in \R^L$ such that $|c|_1\le 1$ we can control
the error of approximation of  $\sum_{\ell} c_\ell \gamma_\ell$ and
$\sum_{\ell} c_\ell \eta_\ell$ with their empirical versions. Namely, we
have on $A$:
\[
%\P\left[
\max_{|c|_1\le 1}
\left|\sum_{\ell}c_\ell(\hat{\eta}_\ell-\eta_\ell)\right|_2\le  \delta_N
%\right]< \varepsilon,
\;\;\mbox{and}\;\;%\P\left[
\max_{|c|_1\le 1}
\left|\sum_{\ell}c_\ell(\hat{\eta}_\ell-\eta_\ell)\right|_2\le  \nu_N.
%\right]< \varepsilon
\]\\[2ex]

{\noindent}Let now  $\hat{U}=[\hat\eta_1,...,\eta_L]$,
$\hat{G}=[\hat\gamma_1,...,\hat\gamma_L ]$.
%We aim to .
%To this end we modify \rf{proj0} as follows: we constrain $c$ to
%and let $H$ be the standard
%hyperoctahedron of $\R^L$: $H=\{c\in \R^L|\;|c\|_1\le 1\}$.
% what allows us to control the
%maximal error of approximation of $Gc$ and $Uc$ with $\hat{G}c$ and $\hat{U}c$.
When substituting $\hat{U}$ and $\hat{G}$ for $U$ and $G$ into  \rf{proj0}
we come to the following minmax problem:
\be
\min_\Pi\max_c\left\{|(I-\Pi)\hat{U}c|^2_2\left|\;\begin{array}{c}\mbox{$\Pi$
is a projectoron an $m$-dimensional}\\ \mbox{subspace of $\R^d$}\\
c\in \R^L,\;|c|_1\le 1,\;|\hat{G}c|_2\le \varrho\end{array}\right.\right\}.
\ee{proj1}
Here we have substituted the constraint $Gc=0$ with the inequality constraint
$|\hat{G}c|_2\le \varrho$ for some $\varrho>0$ in order to keep the optimal
solution $c^*$ to \rf{proj0} feasible for the modified problem \rf{proj1}
(this will be the case with probability at least $1-\varepsilon$ if
$\varrho\ge \nu_N$).

As we will see in a moment, when $c$ runs the
$\nu_N$-neighborhood of intersection $C_N$ of the standard hyperoctahedron
$\{c\in \R^L,\;|c|_1\le 1\}$ with the subspace $\hat{G}c=0$, vectors
$\hat{U}c$ span a close vicinity of the target space $\cI$.

%-------------------------------------------------------------------------
\subsection{Solution by Semidefinite Relaxation}\label{subproj}
%
%Note that for functions $h_{\ell}(x) = f(x,\omega_{\ell})$, $\ell=1,...,L$,
%where $f$ is continuously differentiable function and $\omega_\ell\in \R^d$
%are vectors on the unit sphere, which are the functions used in Section \ref{numerics},
%such that $f(\cdot,\cdot)$ and $\nabla_x f(\cdot,\cdot)$ are bounded,
%\rf{mdp} holds for ``empirical counterparts'' \rf{empir} of $\eta_\ell$
%and $\gamma_\ell$ with $\delta_N=\cO\left(d\sqrt{\min\{d,\ln L\}+\ln\varepsilon^{-1}\over N}\right)$
%(cf., e.g., Lemma 1 in \cite{DJSS}).
%
Note that \rf{proj1} is a hard  optimization problem. Namely, the
candidate maximizers $c_i$ of \rf{proj1} are the extreme points of the set
$C_N=\{c\in \R^L,\,|c|_1\le 1,\, |\hat{G}c|_2\le \nu_N\}$, and there are
$O(L^d)$ of such points. In order to be efficiently solvable, the problem
\rf{proj1} is to be ''reduced'' to a convex-concave saddle-point problem,
which is, to the best of our knowledge, the only class of minmax problems
which can be solved efficiently (cf. \cite{NY}).

Thus the next step is to transform the problem in \rf{proj1}
into a convex-concave minmax problem using the {\em Semidefinite Relaxation} (or
SDP-relaxation) technique (see e.g., \cite[Chapter 4]{BN1}). We obtain the relaxed
version of \rf{proj1} in two steps. First, let us rewrite the objective
function (recall that $I-\Pi$ is also a projector, and thus an idempotent
matrix):
\[
|(I-\Pi)\hat{U}c|^2_2=c^T\hat{U}^T(I-\Pi)^2\hat{U}c=
c^T\hat{U}^T(I-\Pi)\hat{U}c=\tr\left[\hat{U}^T(I-\Pi)\hat{U}X\right],
\]
where the positive semidefinite matrix $X=cc^T$ is the ''new variable''.
The constraints on $c$ can be easily rewritten for $X$:
\begin{enumerate}
 \item the constraint $|c|_1\le 1$ is equivalent to $|X|_1\le 1$
  (we use the notation $|X|_1=\sum_{i,j=1}^L |X_{ij}|$);
 \item because $X$ is positive semidefinite, the constraint $|\hat{G}c|_2\le \varrho$ is equivalent to
  into $\tr[\hat{G}X\hat{G}^T]\le \varrho^2$.
\end{enumerate}
The only ''bad'' constraint on $X$ is  the rank constraint: $\rank X=1$,
and we simply remove it. Now we are done with the variable $c$ and we arrive at
\bse
\min_\Pi\max_X\left\{\tr\left[\hat{U}^T(I-\Pi)\hat{U}X\right]\left|\;
\begin{array}{c}\mbox{$\Pi$ is a projector on an $m$-dimensional }\\ \mbox{subspace of $\R^d$}\\
X\succeq 0,\;|X|_1\le 1,\;\tr[\hat{G}X\hat{G}^T]\le \varrho^2\end{array}\right.\right\}.
\ese
Let us recall that an $m$-dimensional projector $\Pi$ is exactly a
symmetric $d\times d$ matrix of  $\rank \Pi=m$ and $\tr \Pi=m$, with the
eigenvalues $0\le \lambda_i(\Pi)\le 1,\;i=1,...,d$. %Therefore, the
%identification minmax problem from above can be rewritten as
%\be
%\min_\Pi\max_c\left\{\|(I-\Pi)Uc\|^2_2\left|\;\begin{array}{c}0\preceq\Pi\preceq I,\;\tr \Pi=m,\; \rank \Pi=m;\\
%c\in \R^L,\;Gc=0\end{array}\right.\right\}
%\ee{proj0}
%.
%Yet the problem \rf{proj2} is still not convex in $\Pi$. However,
%this can be easily corrected:
Once again we remove the ``difficult'' rank
constraint $\rank \Pi=m$ and finish with
\be
\min_P\max_X\left\{\tr\left[\hat{U}^T(I-P)\hat{U}X\right]\left|\;
\begin{array}{c}0\preceq P\preceq I,\;\tr P=m,\; \\
X\succeq 0,\;|X|_1\le 1,\;\tr[\hat{G}X\hat{G}^T]\le
\varrho^2\end{array}\right.\right\} \ee{proj3} (we write $P\preceq Q$ if
the matrix $Q-P$ is positive semidefinite). There is no reason for an optimal solution
$\hat{P}$ of \rf{proj3} to be a projector matrix. If  an
estimation of  $\Pi^*$ which is itself a projector is needed, one can use instead the
projector $\hat\Pi$ onto the subspace spanned by $m$ principal
eigenvectors of $\hat{P}$.

Note that \rf{proj3} is a linear matrix game with bounded
convex domains of its arguments - positive semidefinite matrices $P\in
\R^{d\times d}$ and $X\in \R^{L\times L}$.

We are about to describe the accuracy of the estimation
$\hat{\Pi}$ of $\Pi^*$. To this end we need an identifiability assumption
on the system $\{h_\ell\}$ of test functions as follows:
\begin{assumption} \label{ass:1}
Suppose that there are vectors $c_{1},...,c_{\bar{m}}$, $m\le \bar{m}\le L$
such that $|c_k|_1\le 1$ and $Gc_k=0$, $k=1,...,\bar{m}$, and
non-negative constants $\mu^{1},\ldots ,\mu^{\bar{m}}$ such that
\be
\Pi^{*} \preceq \sum_{k=1}^{\bar{m}} \mu^{k} Uc_kc_k^TU^T.
\ee{ident}
We denote $\mu^*=\mu^1+\ldots+\mu^{\bar{m}}$.
\end{assumption}
{\noindent} In other words, if Assumption \ref{ass:1} holds, then the true
projector $\Pi^*$ on $\cI$ is $\mu^*\times$ convex combination of rank-one
matrices $Ucc^TU^T$ where $c$ satisfies the constraint $Gc=0$ and
$|c|_1= 1$. % (or, what is the same, the set $\mu^*\conv (Uc)$, where $\conv (Uc)$ is the convex hull of  vectors $Uc$ with
% $Gc=0$ and $|c|_1= 1$, contains an orthogonal basis of $\cI$).
%%%%%%%%%%%%%%%%%%%%%%%%%%%%%%%%%%%%%%
\begin{theorem}\label{the:cond1}
Suppose that the true dimension $m$ of the
subspace $\cI$ is known and that $\varrho\ge \nu_N$ as in \rf{mdp}. Let  $\hat{P}$ be an
optimal solution to \rf{proj3} and let $\hat{\Pi}$ be the projector onto the subspace spanned by $m$
principal eigenvectors of $\hat{P}$. Then with probability $\ge 1-\varepsilon$:
\item{(i)}  for any $c$ such that $|c|_1\le 1$ and $Gc=0$,
\[
|(I-\hat\Pi)Uc|_2\le \sqrt{m+1}((\varrho+\nu_N)\lambda^{-1}_{\min}(\Sigma)+2\delta_N);
\]
\item{(ii)} further, if Assumption \ref{ass:1} holds then
\be
\tr\left[\,(I-\hat{P})\Pi^*\right]\le  \mu^*((\varrho+\nu_N)\lambda^{-1}_{\min}(\Sigma)+2\delta_N)^2,
\ee{id1}
and
%\[
%\tr\left(\,(I-\hat{\Pi})\Pi^*\right)\le {4 \mu^*\delta^2(\lambda^{-1}_{\min}(\Sigma)+1)^2\over 1-4 \mu^*\delta^2(\lambda^{-1}_{\min}(\Sigma)+1)^2},
%\]
%and
\be \begin{array}{rcl}\|\hat{\Pi}-\Pi^*\|_2^2&\le& {2
\mu^*(\lambda^{-1}_{\min}(\Sigma)(\varrho+\nu_N)+2\delta_N)^2} \,\tau, \vspace{0.2cm}\\
\tau&=&(m+1)\wedge (1-
\mu^*(\lambda^{-1}_{\min}(\Sigma)(\varrho+\nu_N)+2\delta_N)^2)^{-1}
\end{array} \ee{id2}
(here
$\|A\|_2=\left(\sum_{i,j}A^2_{ij}\right)^{1/2}=\left(\tr[A^TA]\right)^{1/2}$
is the Frobenius norm of $A$).
\end{theorem}
Note that if we were able to solve the minimax problem in \rf{proj1}, we could expect its solution, let us call it $\tilde\Pi$, to satisfy with high probability
\bse
|(I-\tilde\Pi)Uc|_2\le (\varrho+\nu_N)\lambda^{-1}_{\min}(\Sigma)+2\delta_N
\ese
(cf. the proof of Lemma \ref{obj1} in the appendix). If we compare this bound to that of
the statement (i) of Theorem \ref{the:cond1}, we conclude that the loss of the accuracy resulting from the substitution of \rf{proj1} by its treatable approximation \rf{proj3} is bounded with $\sqrt{m+1}$. In other words, the ``price'' of the SDP-relaxation in our case is $\sqrt{m+1}$ and does not depend on problem dimensions $d$ and $L$. Furthermore, when Assumption 1 holds true, we are able to provide the bound on the accuracy of recovery of projector $\Pi^*$ which is seemingly as good as if we were using instead of $\hat\Pi$ the solution $\tilde\Pi$ of \rf{proj1}.

Suppose now that the test functions $h_\ell(x)=f(x,\omega_\ell)$ are used, with $\omega_l$ on the unit sphere of $\R^d$, that $\varrho=\nu_N$ is chosen, and that Assumption 1 holds with  ``not too large'' $\mu^*$, e.g., $\mu^*\le {\small 1\over 2}(\varrho+\nu_N)\lambda^{-1}_{\min}(\Sigma)+\delta_N$. When substituting the bounds of \rf{dnu} for $\delta_N$ and $\nu_N$ into \rf{id2} we obtain the bound for the accuracy of the estimation $\hat \Pi$ (with probability $1-\epsilon$):
\[
\|\hat{\Pi}-\Pi^*\|_2^2\le C(f)\mu^*N^{-1}\left(\min(d,\ln L)+\ln \epsilon^{-1}\right)
\]
where $C(f)$ depends only on $f$.
This bound claims the root-\( N \) consistency in estimation of the non-Gaussian
subspace with the log-price for relaxation and estimation error.

\subsection{Case of unknown dimension $m$}
The problem \rf{proj3} may be modified to allow the treatment of the case
when the dimension $m$ of the target space is unknown {\em a priori}.
Namely, consider for $\rho\ge0$ the following problem
\be
\min_{P,t} \left\{t\,\left|\;\begin{array}{c}\tr P\le t,\;\max_X\tr\left[\hat{U}^T(I-P)\hat{U}X\right]\le \rho^2,\;0\preceq P\preceq I, \\
X\succeq 0,\;|X|_1\le 1,\;\tr[\hat{G}X\hat{G}^T]\le
\varrho^2\end{array}\right.\right\}
\ee{proj4}
% Note that we can use a
%bisection or a kind of Newton search in $\rho$ (note that the objective of
%\rf{proj4} is obviously convex in $\rho$) to reduce \rf{proj4} to a small
%sequence to feasibility problems, closely related to \rf{proj3}: given
%$t_0$ report (if exists) $P$ such that
%\[
%\max_X\left\{\;\begin{array}{c}\tr\left[\hat{U}^T(I-P)\hat{U}X\right]\le \rho^2,\;0\preceq P\preceq I,\;\tr P\le t_0,\; \\
%X\succeq 0,\;|X|_1\le 1,\;\tr[\hat{G}X\hat{G}^T]\le \varrho^2\end{array}\right\}.
%\]
%{\noindent}In other words, we can easily solve \rf{proj4} if for a given
%$m$ we are able to find an optimal solution to \rf{proj3}.

The problem \rf{proj4} is closely related to the
$\ell_1$-recovery estimator of sparse signals (see, e.g., the tutorial
\cite{candes} and the references therein) and the trace minimization
heuristics widely used in the Sparse Principal Component Analysis (SPCA)
(cf. \cite{AE07,ABE08}). As we will see in an instant, when the parameter
$\rho$ of the problem is ''properly chosen'', the optimal solution
$\hat{P}$ of \rf{proj4} possesses essentially the same properties as that
of the problem \rf{proj3}.

A result analogous to that in Theorem \ref{the:cond1} holds:
\begin{theorem}\label{the:cond2}
Let $\hat{P},
\;\hat X$ and $\hat t=\tr \hat P$ be an optimal solution to \rf{proj4} (note that \rf{proj4} is clearly solvable), $\hat m=\rfloor \hat t\lfloor$,\footnote{Here $\rfloor a\lfloor$ is the smallest integer $\ge a$.} and let
$\hat{\Pi}$ be the projector onto the subspace spanned by $\hat m$
principal eigenvectors of $\hat{P}$. Suppose that  $\varrho\ge \nu_N$ as in
\rf{mdp} and that \be \rho\ge
\lambda^{-1}_{\min}(\Sigma)(\varrho+\nu_N)+\delta_N. \ee{rhobound}
Then with probability
at least $1-\varepsilon$:
\item{(i)}
\[
\hat{t}\le m\;\;\mbox{and}\;\;|(I-\hat\Pi)Uc|_2\le \sqrt{m+1}(\rho+2\delta_N);
\]
\item{(ii)} furthermore, if Assumption \ref{ass:1} hold then
\[
\tr\left[\,(I-\hat{P})\Pi^*\right]\le  \mu^*(\rho+\delta_N)^2,
\]and\
\be
\|\hat{\Pi}-\Pi^*\|_2^2\le {2\mu^*(\rho+\delta_N)^2}\,\left[(m+1)\wedge (1- \mu^*(\rho+\delta_N)^2)^{-1}\right]
\ee{id22}
(here $\|A\|_2=\left(\sum_{i,j} A_{ij}^2\right)^{1/2}$  is the Frobenius norm of $A$).
\end{theorem}
The proof of the theorems is postponed until the
appendix.

The estimation procedure based on solving \rf{proj4} allows to
infer the target subspace $\cI$ without {\em a priori} knowledge of its
dimension $m$. When the constraint parameter $\rho$ is close to the
right-hand side of \rf{rhobound}, the accuracy of the estimation will be
close to that, obtained in the situation when dimension $m$ is known.
However, the accuracy of the estimation heavily depends on the precision
of the available (lower) bound for $\lambda_{\min}(\Sigma)$. In the
high-dimensional situation this information is hard to acquire, and the
necessity to compute this quantity may be considered as a serious drawback
of the proposed procedure.

\section{Solving the saddle-point problem \rf{proj3}}
We start with the following simple observation: by using bisection or Newton search in $\rho$ (note that the objective of \rf{proj4} is obviously convex in $\rho^2$) we can reduce \rf{proj4} to a small sequence to feasibility problems, closely related to \rf{proj3}: given $t_0$ report, if exists, $P$ such that
\[
\max_X\left\{\;\begin{array}{c}\tr\left[\hat{U}^T(I-P)\hat{U}X\right]\le \rho^2,\;0\preceq P\preceq I,\;\tr P\le t_0,\; \\
X\succeq 0,\;|X|_1\le 1,\;\tr[\hat{G}X\hat{G}^T]\le \varrho^2\end{array}\right\}.
\]
In other words, we can easily solve \rf{proj4} if for a given $m$ we are  able to find an optimal solution to \rf{proj3}.
Therefore, in the sequel we concentrate on the optimization technique for solving \rf{proj3}.
\subsection{Dual extrapolation algorithm}
In what follows we discuss the dual extrapolation algorithm of \cite{N3} for solving a version of \rf{proj3} in which, with a certain abuse, we substitute the inequality constraint $\tr \hat G X\hat G^T\le \varrho^2$ with the equality constraint $ \tr [\hat G X\hat G^T]=0$. This way we come to the problem:
\be
\min_{P\in {\cal P}}\max_{X\in \cal X} \tr\left[\hat{U}^T(I-P)\hat{U}X\right]
\ee{proj5}
where
\begin{eqnarray*}
{\cal X}=\{X\in {\bold S}^L, \;X\succeq 0,\;|X|_1\le 1,\;\tr [\hat G^T\hat G X]=0\}
\end{eqnarray*}
(here ${\bold S}^L$ stands for the space of $L\times L$ symmetric matrices) and
\begin{eqnarray*}
{\cal P}=\{P\in {\bold S}^d,\;0\preceq P \preceq I,\; \tr[ P]\le m\}.
\end{eqnarray*}
Observe first that \rf{proj5} is a matrix game over two convex subsets (of the cone) of positive semidefinite matrices.
If we use a large number of test functions, say $L^2\sim 10^6$, the size of the variable $X$ rules
out the possibility of using the interior-point methods. The methodology which appears
to be adequate in this case is that behind dual extrapolation methods,
recently introduced in \cite{N1,N12,N2,N3}. The algorithm we use
belongs to the family of subgradient descent-ascent
methods for solving convex-concave games. Though the rate of convergence of
such methods is slow --- their precision is only \(\cO( 1/k )\),
where $k$ is the iteration count, their iteration is relatively cheap, what makes the methods
of this type appropriate in the case of high-dimensional
problems when the high accuracy is not required.

%We are to write the dual extrapolation method of \cite{N3}  for the problem \rf{proj5}.
We start with the general dual extrapolation scheme of \cite{N3} for linear matrix games. Let $\eE^n$ and $\eE^m$ be two Euclidean spaces of dimension $n$ and $m$ respectively, and let $\cA\subset \eE^n$ and $ \cB\subset \eE^m$ be closed and convex sets. We consider the problem
\be
\min_{x\in \cA}\max_{y\in \cB}\langle x,Ay\rangle+\langle a,x\rangle+\langle b,y\rangle.
\ee{sp}
Let $\|\cdot\|_x$ and $\|\cdot\|_y$ be some norms on $\eE^n$ and $\eE^m$ respectively. We say that $d_x$ (resp., $d_y$) is a distance-generating function of $\cA$ (resp., of $\cB$) if $d_x$ (resp., $d_y$) is strongly convex modulus $\alpha_x$ (resp., $\alpha_y$)  and differentiable on $\cA$ (resp., on $\cB$).\footnote{Recall that a (sub-)differentiable on $\cF$ function $f$ is called strongly convex on $\cF$ with respect to the norm $\|\cdot\|$ of modulus $\alpha$ if  $\langle f'(x)-f'(y), x-y\rangle \ge \alpha \|x-y\|^2$ for all $x,y\in \cF$.} Let for
 $z=(x,y)$  $d(z)=d_x(x)+d_y(y)$ (note that $d$ is  differentiable and strongly convex on $\cA\times \cB$ with respect to the norm, defined on $\cA\times \cB$ according to, e.g. $\|z\|=\|x\|_x+\|y\|_y$). We define the {\em prox-function} $V$ of $\cA\times \cB$ as follows: for  $z_0=(x_0,y_0)$ and $z=(x,y)$ in $\cA\times \cB$ we set
\be
V(z_0,z)\eqdef d(z)-d(z_0)-\langle \nabla d(z_0),z-z_0\rangle.
\ee{proxf}
Next, for $s=(s_x,s_y)$ we define the {\em prox-tranform} $T(z_0,s)$ of $s$:
\be
T(z_0,s)\eqdef\arg\min_{z\in \cA\times \cB}[\langle s,\;z-z_0\rangle- V(z_0,z)].
\ee{proxt}
Let us denote $F(z)=(-A^Ty-a,Ax+b)$ the vector field of descend-ascend directions of \rf{sp} at $z=(x,y)$ and let $\bar{z}$ be the minimizer of $d$ over $\cA\times \cB$. Given vectors $z_k,z^+_k\in \cA\times \cB$ and $s_k\in E^*$ at the $k$-th iteration, we define the update $z_{k+1},\,z^+_{k+1}$ and $s_{k+1}$ according to
\bse
z_{k+1}&=&T(\bar{z},s_k),\\
z^+_{k+1}&=&T(z_{k+1},\lambda_k F(z_{k+1})),\\
s_{k+1}&=&s_k+\lambda_k F(z^+_{k+1}),
\ese
where $\lambda_k>0$ is the current stepsize. Finally, the current approximate solution $\hat{z}_{k+1}$ is defined with
\[
\hat z_{k+1}={1\over k+1}\sum_{i=1}^{k+1} z^+_{i}.
\]
The key element of the above construction is the choice of the {\em distance-generaing function} $d$ in the definition of the prox-function. It should satisfy two requirements:
\begin{itemize}
\item let $D$ be the variation of $V$ over $\cA\times \cB$ and let $\alpha$ be the parameter of strong convexity of $V$ with respect to $\|\cdot\|$.
    The complexity of the algorithm is proportional to $D/\alpha$, so this ratio should be as small as possible;
\item one should be able to compute efficiently the solution to the auxiliary problem \rf{proxt} which is to be solved twice at each iteration of the algorithm.
\end{itemize}
Note that the prox-transform preserve the additive structure of the distance-generating function. Thus, in order to compute the prox-transform on the feasible domain $\cP\times \cX$ of \rf{proj5}  we need to compute its ``$P$ and $X$ components'' -- the corresponding prox-transforms on $\cP$ and $cX$. There are several evident choices of the prox-functions $d_P$ and $d_X$ of the domains $\cP$ and $\cX$ of \rf{proj5} which satisfy the first requirement above and allow to attain the optimal value $O(\sqrt{m\ln d\ln L})$ of the ratio ${D/ \alpha}$ for the prox-function $V$ of \rf{proj5}.
However, for such distance-generating functions there is no known way to compute  efficiently the $X$-component of the prox-transform $T$ in \rf{proxt} for the set $\cal X$.
This is why in order to admit an efficient solution  the problem \rf{proj5} is to be modified one more time.
\subsection{Modified problem}
We act as follows: first we eliminate the linear equality constraint which, taken along
with $X\succeq 0$, says that $X = Q^TZQ$ with $Z \succeq 0$ and certain $Q$; assuming
that the $d$ rows of $\hat G$ are linearly independent, we can choose $Q$ as an
appropriate $(L-d)\times L$ matrix satisfying $QQ^T = I$ (the orthogonal basis of the
kernel of $\hat G$). Note that from the constraints on $X$ it follows that $\tr [X]\le
1$, whence \begin{eqnarray*}\tr [Q^TZQ] = \tr [ZQQ^T] = \tr [Z]\le 1.\end{eqnarray*}
Thus, although there are additional constraints on Z as well, $Z$ belongs to the {\em
standard spectahedron} \bse{\cal Z} = \{Z\in{\bold S}^{L-d},\; Z\succeq 0,\,\tr [Z]\le
1\}.
 \ese
Now can rewrite our problem equivalently as follows:
\be
\min_{P\in \cP}
\max_{%\begin{array}{c}
Z\in \cZ,\,|Q^TZQ|_1\le 1
%\end{array}
}
\tr[\hat U^T(I-P) \hat U(Q^TZQ)].
\ee{arikoo}
Let, further,
\bse
\cW = \{W\in {\bold S}^{L}, \,\|W\|_2\le 1\},\;\;\mbox{and}\;\; \cY = \{Y\in  {\bold S}^{L},\,|Y|_1\le 1\}.
\ese
We claim that the problem \rf{arikoo} can be reduced to the saddle point
problem
\be
\min_{(P,W)\in \cP\times \cW}
\max_{(Z,Y)\in \cZ\times \cY}
\underbrace{\left\{
\tr[\hat U^T (I-P)\hat UY]
+\lambda\, \tr[W(Q^TZQ-Y)]\right\}
}_{F(P,W;\,Z,Y)}.
\ee{2}
provided that $\lambda$ is not too small.

Now, ``can be reduced to'' means exactly the following:
\begin{proposition}\label{arik100}
Suppose that $ \lambda>L|\hat U|_2^2$, where $|U|_2$ is the maximal Euclidean norm of columns of $U$. Let $(\hat P, \hat W;\,\hat Z,\hat Y)$ be a feasible solution $\e$-solution to \rf{2}, that is
\bse
(\hat P, \hat W;\,\hat Z,\hat Y)\in  (\cP, \cW; \cZ, \cY),\;\; \mbox{and}\;\;\bar F(\hat P, \hat W)-\underline{F}(\hat Z,\hat Y)\le \e
\ese
where
\bse
\bar F(P,W) = \max_{(Z,Y)\in \cZ\times \cY}
F(P,W;\,Z, Y),\;\;\; \underline{F}(Z,Y) = \min_{(P,W)\in \cP\times \cW}
F(P,W;\,Z, Y).
\ese
Then setting
\bse
\tilde{Z} =\left\{\begin{array}{ll}
\hat{Z},&\mbox{if}\; |Q^TZQ|_1\le 1,\\|Q^TZQ|_1^{-1}\hat{Z}&\mbox{otherwise},\end{array}\right.
\ese
the pair $(\hat P, \tilde{Z})$ is a feasible $\e$-solution to \rf{arikoo}. Specifically, we have $(\hat P, \tilde Z)\in  \cP \times \cZ$ with $|Q^T\tilde Z Q|_1\le 1$, and
\bse
\bar G(\hat P)- \underline{G}(\tilde Z)\le \e,
\ese
where
\bse
\bar G(P) &=& \max_{Z\in \cZ,\,|Q^TZQ|_1\le 1}
\tr[\hat U^T (I-P)\hat U Q^TZQ]; \;\;\underline{G}(Z) \\
&=& \min_{P\in \cP}
\tr [\hat U^T (I-P)\hat U Q^TZQ].
\ese
\end{proposition}
The proof of the proposition is given in the appendix \ref{sec:arik}.

Note that feasible domains of \rf{2} admit
evident distance-generating functions. We provide the detailed computation of the corresponding prox-transforms in the appendix \ref{sec:prox}.

%-------------------------------------------------------------------------
\section{Numerical Experiments}\label{numerics}
%-------------------------------------------------------------------------
In this section we compare the numerical performance of the presented
approach, which we refer to as SNGCA(SDP) with other statistical methods of dimension reduction on
the simulated data.

\subsection{Structural adaptation algorithm}
We start with some implementation details of the
estimation procedure.
We use the choice of the test functions $h_\ell(x)=f(x,\omega_\ell)$ for the SNGCA algorithm as follows:
\bse
    f(x,\omega)
    &=&
    \tanh(\omega^{T}x) e^{-\alpha\|x\|^{2}_{2}/2},
\ese where  $\omega_{\ell}$, $l=1,...,L$ are  unit vectors in
$\R^{d}$.

We implement here a multi-stage variant of the SNGCA
(cf \cite{DJSS}). At the first stage of the SNGCA(SDP) algorithm we assume that
the directions $\omega_{\ell}$ are drawn randomly from the unit sphere of
$\R^{d}$. At each of the following stages  we use the current estimation
of the target subspace to ``improve'' the choice of directions
$\omega_\ell$ as follows: we draw a fixed fraction of $\omega$'s from the
estimated subspace and draw randomly over the unit square the remaining
$\omega$'s. The simulation results below are present for the estimation
procedure  with three stages. The size of the set of test function is set to $L=10\,d$, and the target accuracy of solving the problem \rf{proj3} is set to $1e-4$.
\par
We can summarize the SNGCA(SDP) algorithm as follows:

\begin{algorithm2e}[H]\label{alg:full}
%\SetLine \dontprintsemicolon \KwData{$\{X_i\}_{i=1}^N$,$L$,$m$}
%\KwResult{$\widehat{\cI}$}%
%\BlankLine
{\tt \% Initialization:}
\par \noindent The data \( (X_{i})_{i=1}^{N} \) are
re-centered. Let \(\sigma=(\sigma_1,\ldots\sigma_{d}) \) be the standard
deviations of the components of \( X_{i} \). We denote  \(
Y_{i}=\diag(\sigma^{-1})X_{i} \) the standardized data.

Set the current estimator $\hat \Pi_0=I_d$.
\BlankLine
{\tt \% Main iteration loop:} \BlankLine \For{i=1 \KwTo I}{
  \BlankLine
  Sample  a fraction of \( \omega^{(i)}\)'s  from the normal distribution
  $N(0,\hat \Pi_{i-1})$ (zero mean, with covariance matrix $\hat \Pi_{i-1}$),
  sample the remaining \(\omega^{(i)}\)'s from $N(0,I_d)$, then normalize to the unit length;
  \BlankLine
  {\tt \% Compute estimations of $\eta_\ell$ and $\gamma_\ell$}

     \For{$\ell$=1 \KwTo L}
     {
       $\widehat{\eta}_{\ell}^{(i)}=\frac{1}{N}\sum_{j=1}^N\nabla h_{\omega_{\ell}^{(i)}}(Y_j)$\;
       $\widehat{\gamma}_{\ell}^{(i)}=\frac{1}{N}\sum_{j=1}^N Y_j h_{\omega_{l}^{(i)}}(Y_j)$\;
     }\BlankLine
     Solve the corresponding problem \rf{proj3} and update the estimation $\hat \Pi_{i}$;
   \BlankLine
}
 \caption{SNGCA (SDP)}
\end{algorithm2e}

\subsection{Experiment description}
\begin{figure}[H]
  \centering
  \begin{tabular}{@{} ccc@{}}
    \includegraphics[width=0.255\textwidth]{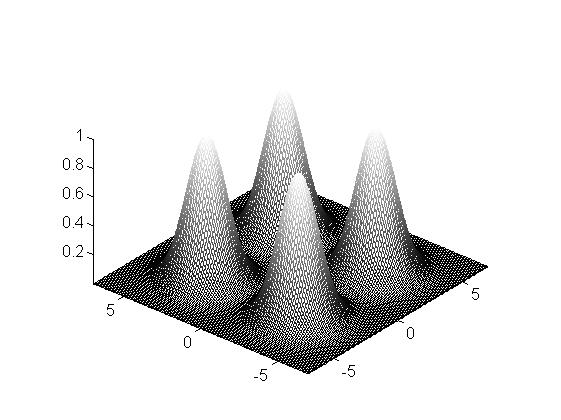} &
    \includegraphics[width=0.255\textwidth]{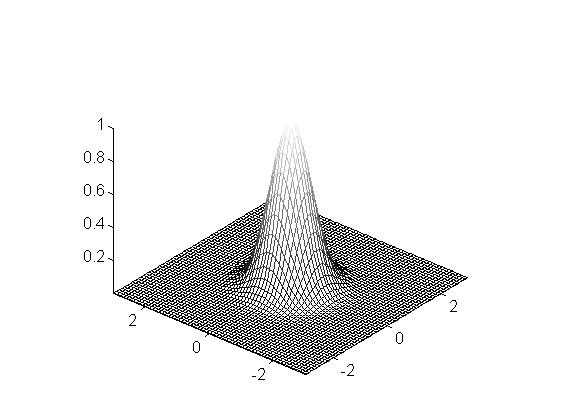} &
    \includegraphics[width=0.255\textwidth]{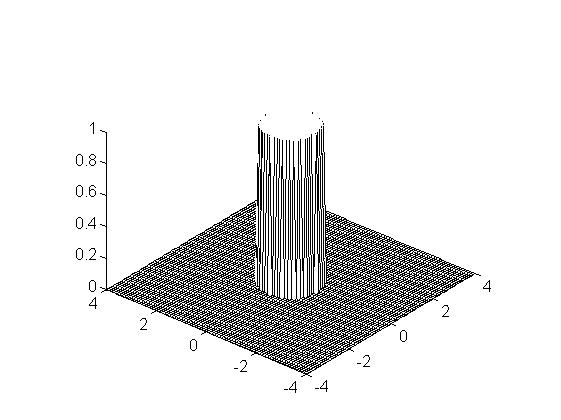} \\
        (A) & (B) & (C)\\
    \includegraphics[width=0.255\textwidth]{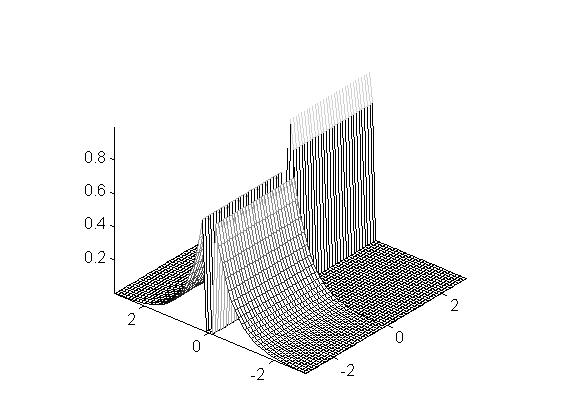} &
    \includegraphics[width=0.255\textwidth]{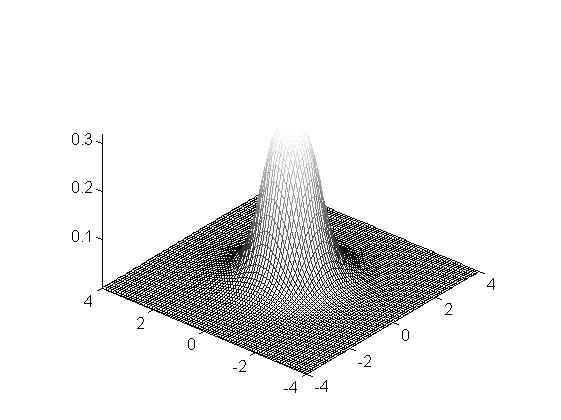} & \\
        (D) & (E) & \\
  \end{tabular}
  \caption{(A) independent Gaussian mixtures, (B) isotropic
    super-Gaussian, (C) isotropic uniform and (D) dependent $1$d Laplacian with additive $1
   $d uniform, (E) isotropic sub-Gaussian  \label{fig:density}
%   {\color{red} Why using different names here and in the description of A-E?}
  }
\end{figure}

Each simulated data set $X^N=[X_1,...,X_N]$ of size $N=1000$ represents
$N$ i.i.d. realizations of a random vectors $X$ of dimension $d$.  Each
simulation is repeated $100$ times and we report the average over $100$
simulations Frobenius norm of the error of estimation of the projection on
the target space. In the examples below only $m=2$ components of $X$ are
non-Gaussian with unit variance, other $d-2$ components of $X$ are
independent standard normal r.v.. The densities of the non-Gaussian
components are chosen as follows:
\begin{description}
\item[(A) Gaussian mixture:] $2$-dimensional independent Gaussian mixtures
with density of each component given by
$0.5\;\phi_{-3,1}(x)+0.5\;\phi_{3,1}(x)$.
\item[(B) Dependent super-Gaussian:] $2$-dimensional isotropic distribution with
density proportional to $\exp(-\|x\|)$.
\item[(C) Dependent sub-Gaussian:] $2$-dimensional isotropic uniform with
constant positive density for $\|x\|_{2}\le 1$ and $0$ otherwise.
\item[(D) Dependent super- and sub-Gaussian:] a component of $X$, say $X_1$, follows the Laplace
distribution $\cL(1)$ and the other is  a
dependent uniform $ \cU(c,c+1) $, where $c=0$ for $|X_1|\le \ln \,2$
and $c=-1$ otherwise.
\item[(E) Dependent sub-Gaussian:] $2$-dimensional isotropic Cauchy distribution
with density proportional to $\lambda(\lambda^2-x^2)^{-1}$ where
$\lambda=1$.
%{\color{red} WHAT IS THIS? Two independent Cauchy or what? }
\end{description}
We provide the $2$-d plots of the densities of the non-Gaussian components
on Figure \ref{fig:density}.

\begin{figure}[H]
    \begin{center}
      \includegraphics[width=13cm]{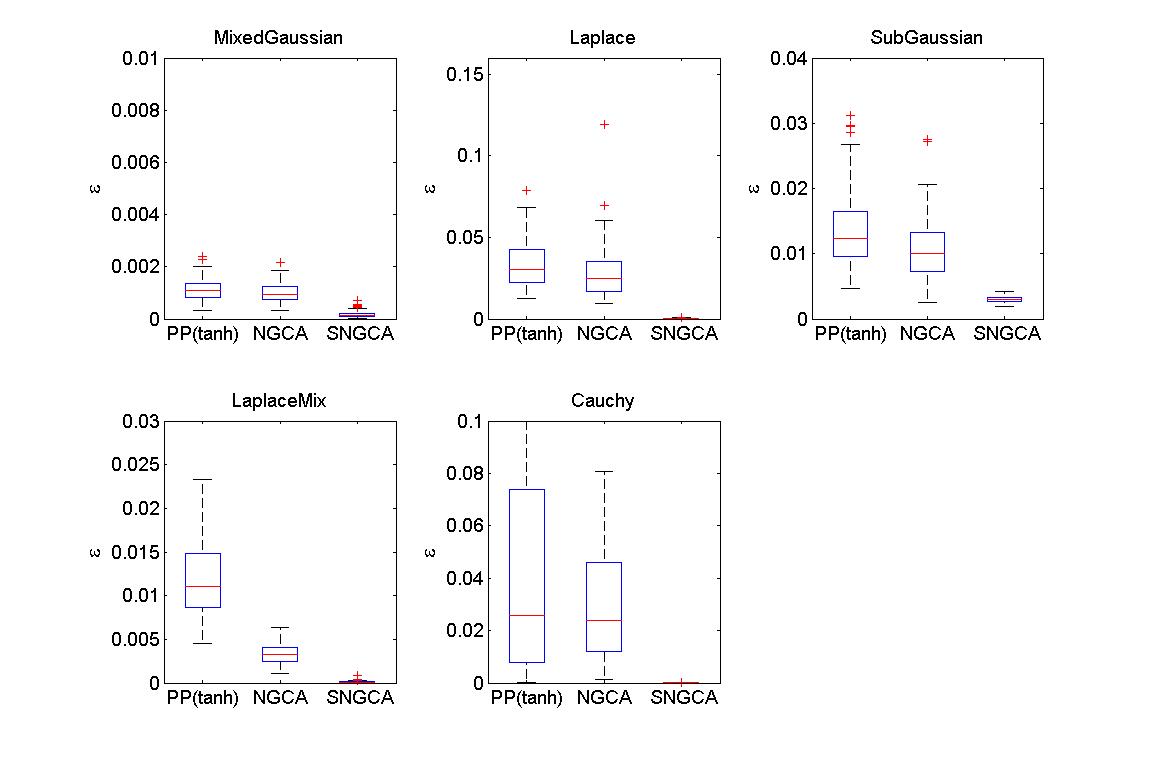}
      \caption{Comparison of PP, NGCA and SNGCA(SDP)}
      \label{fig:AllIterativeMethodsBox}
    \end{center}

\end{figure}

We start with comparing the presented algorithm with Projection Pursuit (PP) method \cite{Hyvaer99} and the NGCA for $d=10$. The results are
 presented on Figure
\ref{fig:AllIterativeMethodsBox} (the corresponding results for PP and NGCA has been already reported in \cite{DJSS} and
\cite{BKSS}).

{\noindent} Since the minimization procedure of PP tends to be trapped  in a local minimum,
in each of the 100 simulations, the PP algorithm  is restarted 10 times  with random
starting points. The best result is reported for each PP-simulation. We
observe that SNGCA(SDP) outperforms  NGCA and PP in all tests.
%Compared to the results obtained from the convex
%projection approach to SNGCA called {\em SNGCA(QCP)} presented in
%\cite{Died07} this is an improvement.\\[2ex]

\begin{figure}[H]
  \centering
  \includegraphics[width=13cm]{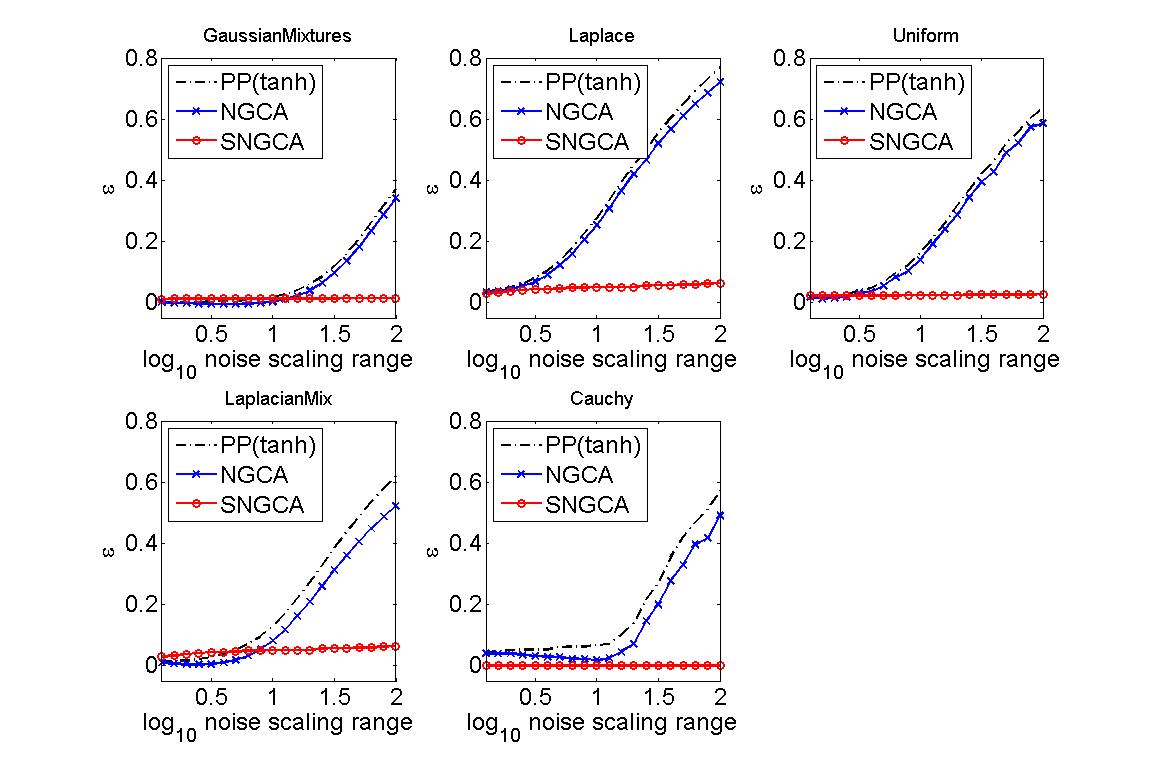}
  \caption{estimation error with respect to the standard deviation of Gaussian components
    following a geometrical progression on $[10^{-r},10^r]$ where $r$ is the parameter
    on the abscissa}
  \label{fig:Noise}
\end{figure}

In the next simulation we study the dependence of the accuracy
of the SNGCA(SDP) on the noise level and compare it to the corresponding data
for PP and NGCA. We present on Figure \ref{fig:Noise} the results of
experiments when the non-Gaussian coordinates have unit variance,
but the standard deviation of the components of the $8$-dimensional
Gaussian distribution follows the geometrical progression
$10^{-r},10^{-r+2r/7},\ldots,10^r$ where $r=1,\ldots,8$.
The conditioning of the covariance matrix
heavily influences the estimation error of PP(tanh) and NGCA, but not that
of SNGCA(SDP). The latter method appears to be insensitive  to the
differences in the noise variance along different direction in all test cases.
% Compared
%to the results of SNGCA(QCP), reported in \cite{Died07}, this is
%also an improvement based on the better use of the information about $\cI$
%obtained from the sampling.\\

\begin{figure}[H]
  \centering
  \includegraphics[width=13cm]{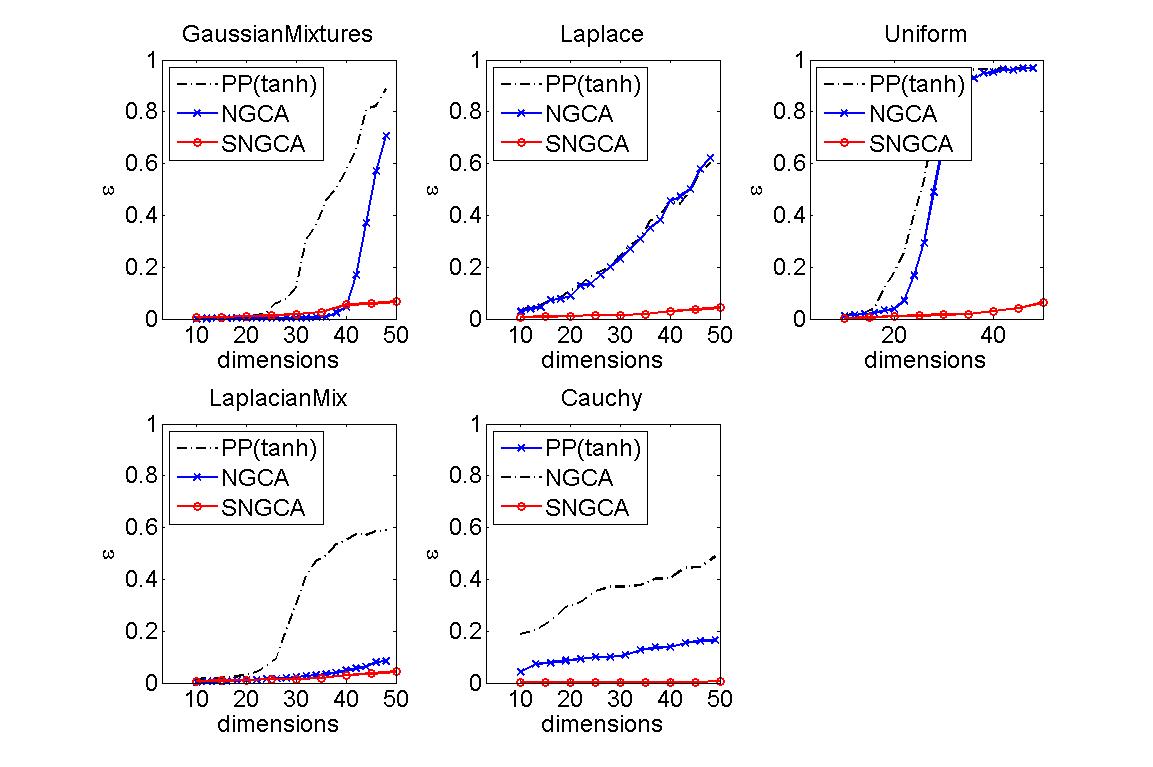}
  \caption{mean-square estimation error vs problem dimension $d$}
  \label{fig:Dim}
\end{figure}

Next we compare the behavior of SNGCA(SDP), PP and NGCA as the
dimension of the  Gaussian component increases. On Figure
\ref{fig:Dim} we plot the mean error of estimation
against the problem dimension $d$.

%{}In figure \ref{fig:dim-increase} we repeat these simulations
%from above in order to demonstrate the specific improvement made by
%SNGCA(SDP) compared to SNGCA(QCP) wrt. to the statistical estimation
%error.
%
%\begin{figure}[H]
%\begin{center}
%  \begin{tabular}{@{} ccc@{}}
%    \includegraphics[width=0.29\textwidth]{MixedGaussianDim.pdf}&
%    \includegraphics[width=0.29\textwidth]{LaplaceDim.pdf} &
%    \includegraphics[width=0.29\textwidth]{UniformDim.pdf} \\
%        (A) & (B) & (C)\\
%        \includegraphics[width=0.29\textwidth]{LaplaceMixDim.pdf} &
%        \includegraphics[width=0.29\textwidth]{CauchyDim.pdf} \\
%        (D) & (E) & \\
%  \end{tabular}
%    \caption{results wrt. $\cE(\widehat{\cI},\cI)$ with increasing number of gaussian components.
%    \label{fig:dim-increase}}
%\end{center}
%\end{figure}

For PP and NGCA methods we observe that the estimation
becomes meaningless (the estimation error explodes) already for $d=30-40$
for the models (A), (C) and for $d=20-30$ of the model  (D). In the case
of the models (B) and (E) we observe the progressive increase of the error
for methods PP and NGCA. The proposed method SNGCA(SDP) behaves robustly
with respect to the increasing dimension of the Gaussian component for all
test models.

%Finally, to illustrate the ``statistical limits'' of
%SNGCA(SDP),  we present on Figure \ref{fig:limits} the results of
%simulations for high-dimensional Gaussian component. In the underlying experiments the number $L$ of test functions \( h_\ell\) is proportional to the
%dimension $d$ {\color{red} HOW MANY EXACTLY?}, and the Gaussian component follows $d-2$-dimensional  standard normal  distribution.
%
%\begin{figure}[H]
%  \centering
%  \includegraphics[width=6in]{IterativeSNGCA-SDPDim}
%    \caption{Estimation error vs  dimension of the problem dimension.
%    \label{fig:limits}}
%\end{figure}
%-------------------------------------------------------------------------
\subsection{Application to Geometric Analysis of Metastability}
%------------------------------------------------------------------------
Some biologically active molecules exhibit different large geometric
structures at the scale much larger than the diameter of the
atoms. If there are more than one
such structures with the life span much larger that the time
scale of the local atomic vibrations, the structure is
called metastable conformation~\cite{a4-SH}.
In other words, metastable
conformations of biomolecules can be seen as connected
subsets of state-space.  When compared to the
fluctuations within each conformation, the
transitions between different conformations of a molecule are rare
statistical events. Such multi-scale dynamic behavior of biomolecules stem
from a decomposition of the free energy landscape into particulary
deep wells each containing many local
minima~\cite{Pillardy95,Frauenfelder00}. Such wells represent
different almost invariant geometrical large scale
structures~\cite{Amadei1993}.
The macroscopic dynamics is
assumed to be a Markov jump process, hopping between the
metastable sets of the state space while the microscopic dynamics
within these sets mixes on much shorter time
scales~\cite{Schutte2008}.
Since the shape of the energy landscape and the invariant density of the Markov process  are unknown, the ``essential degrees
of freedom'', in which the rare conformational changes occur, are of importance.

We will now illustrate that SNGCA(SDP) is able to detect a
multimodal component of the data density as a special case of
non-Gaussian subspace in high-dimensional data obtained from molecular
dynamics simulation of oligopeptides.

\paragraph{Clustering of 8-alanine} The first example is a times series, generated by an
equilibrium molecular dynamics simulation of 8-alanine.
% (courtesy of F. Noe).
%{\color{red} a reference here? Private communication?}.
We only consider the backbone dihedral angles in order to
determine different conformations.
\par
The $14$-dimensional time series consists of the cyclic
data set of all backbone torsion angles. The simulation using {\tt CHARMM} was done at $T=300K$ with
implicit water by means of the solvent model
{\tt ACE2}~\cite{Schaefer96}. A symplectic Verlet integrator with integration
step of $1 fs$ was used; the total
trajectory length was $4 \mu s$ and every
$\tau=50fs$ a set of coordinates was recorded.
\par
The dimension reduction reported in the next figure
was obtained using SNGCA(SDP) with for a given dimension $m=5$ of the
target space containing the multimodal component.

\begin{figure}[H]
    \begin{center}
      \includegraphics[height=7.5cm]{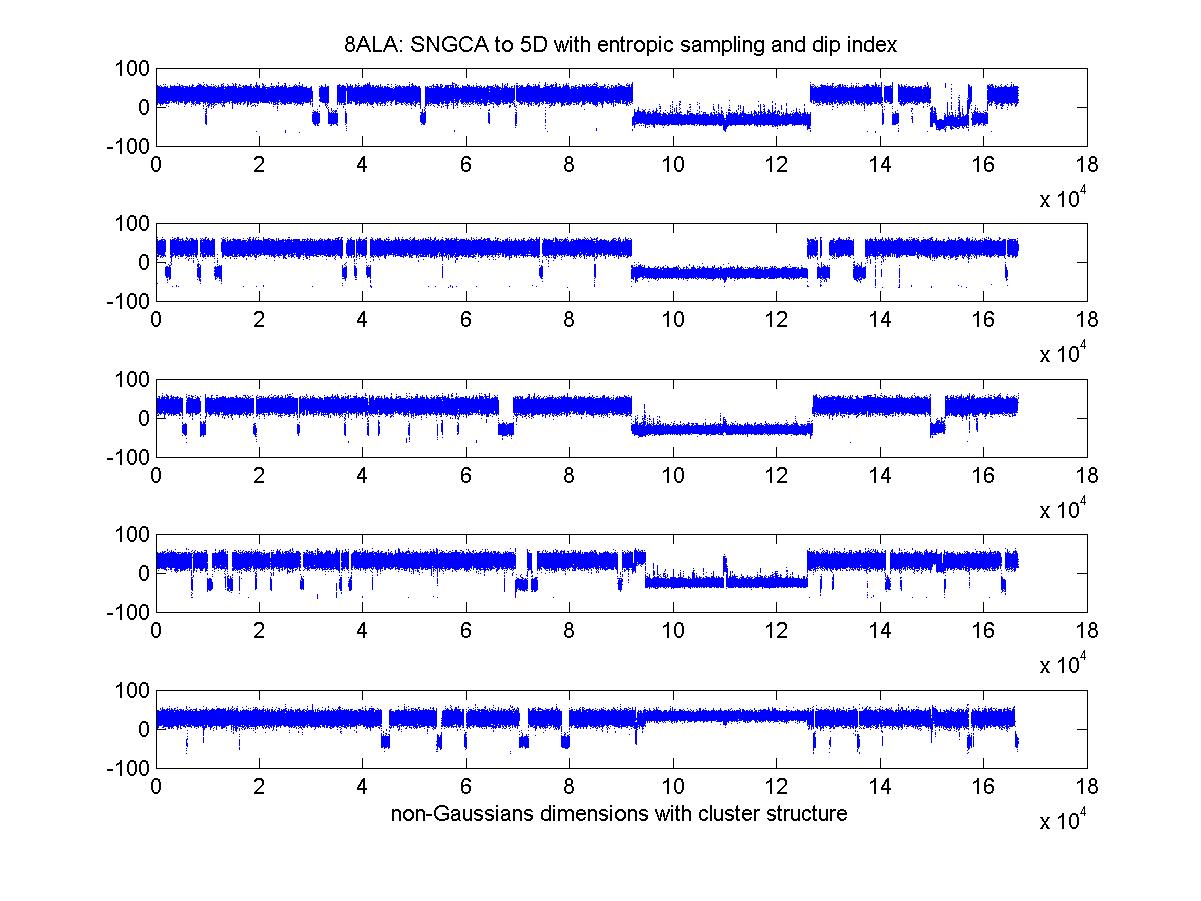}
      \caption{low dimensional multimodal component of 8-alanine
%      {\color{red} Do you understand this figure?}
      \label{fig:8ALA-target}}
    \end{center}
\end{figure}
A concentration of the
clustered data in the target space of SNGCA may be clearly observed. In comparison, the
complement of the target space is almost completely filled with
Gaussian noise.
%\begin{comment}
%\begin{figure}[H]
%    \begin{center}
%      \includegraphics[height=7.5cm]{8ala-SNGCA-complement4}
%      \caption{Gaussian noise in the complement of the SNGCA
%      target space {\color{red} What do you see here?}
%      \label{fig:8ALA-complement1}}
%    \end{center}
%\end{figure}
%\end{comment}
\begin{figure}[H]
    \begin{center}
      \includegraphics[height=7.5cm]{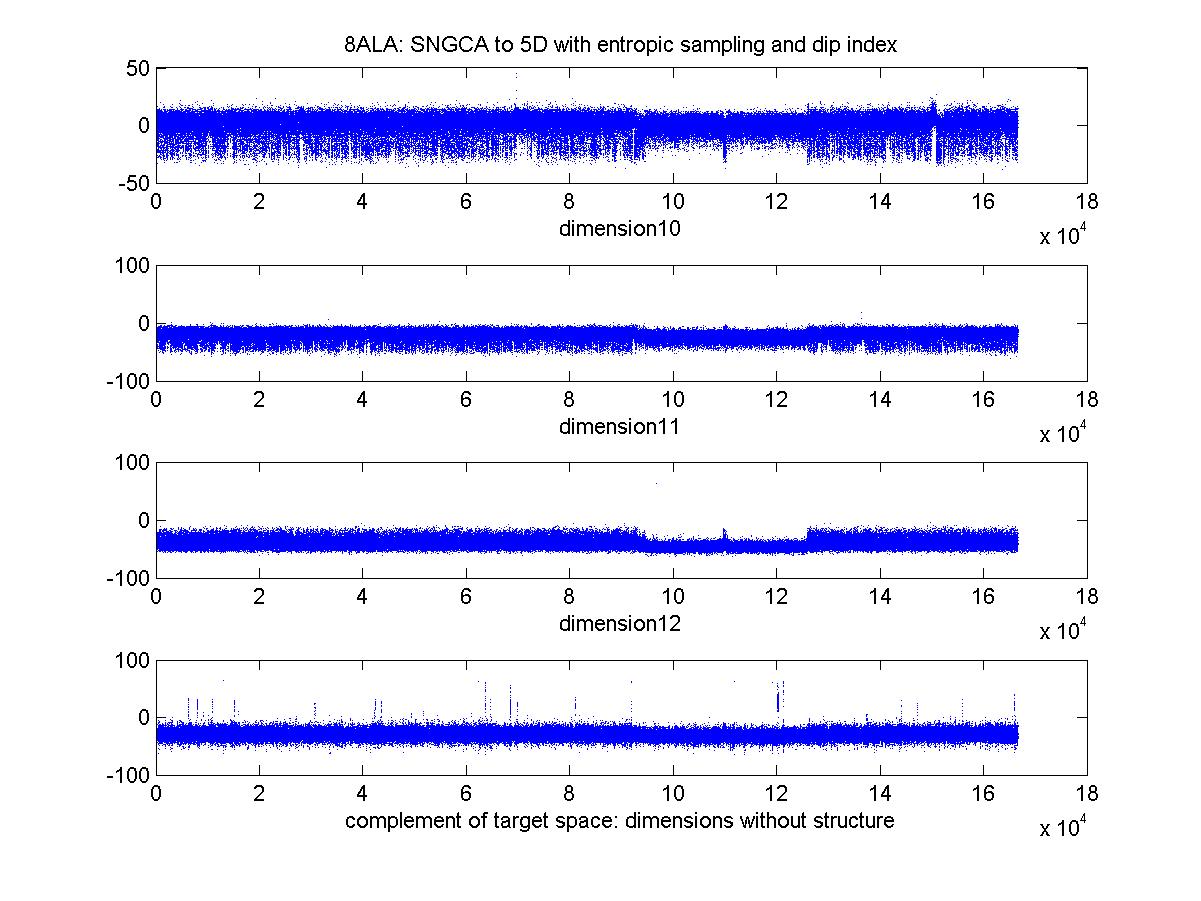}
      \caption{Gaussian noise in the complement of the SNGCA
      target space.
%      {\color{red} What do you see here? BTW, as far as I recall, 5+8=13 and not 14!}
      \label{fig:8ALA-complement2}}
    \end{center}
\end{figure}

\paragraph{Clustering of a 3-peptide molecule} In the next example we investigate
Phenylalanyl-Glycyl-Glycine Tripeptide, which is
assumed to realize all of the most important folding
mechanisms of polypeptides~\cite{RVVVHACV}. The simulation is done
using {\tt GROMACS} at $T=300K$ with implicit water. An integration step
of a symplectic Verlet integrator is set to $2 fs$, and every
$\tau=50fs$ a set of 31 diedre angles was recorded. As in the previous experience, the dimension of the target space is set to $m=5$
\par
\begin{figure}[H]
    \begin{center}
      \includegraphics[height=7.5cm]{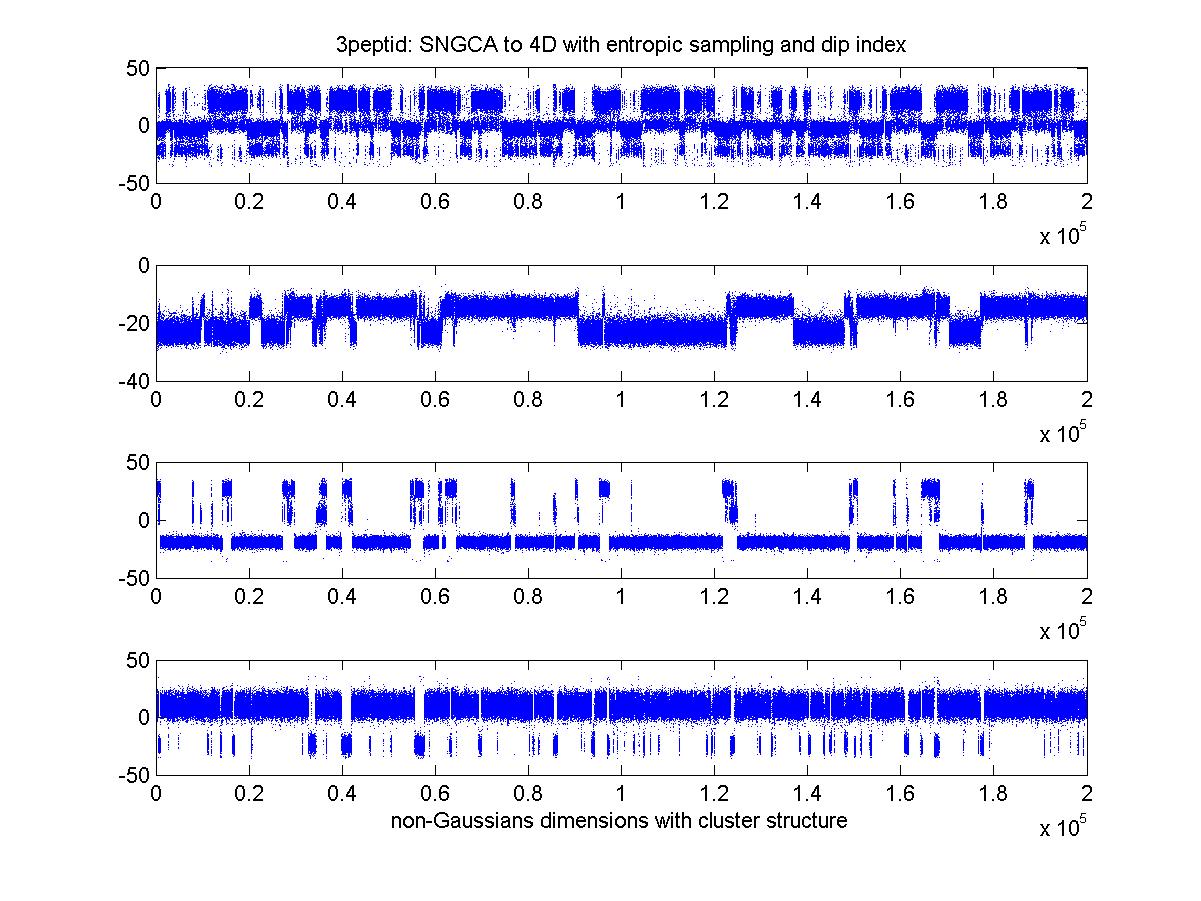}
      \caption{low dimensional multimodal component of 3-peptide
      \label{fig:3pep-target}}
    \end{center}
\end{figure}

Figure \ref{fig:3pep-target} shows that the
clustered data can be primarily found in the target space of
SNGCA(SDP).

\section{Conclusions}
We have studied a new procedure of non-Gaussian component analysis. The suggested method, same as the techniques proposed in \cite{BKSS,DJSS}, has two stages: on the first stage certain linear functionals of unknown distribution are computed, then this information is used to recover the non-Gaussian subspace. The novelty of the proposed approach resides in the new method of non-Gaussian subspace identification, based upon semidefinite relaxation. The new procedure allows to overcome the main drawbacks of the previous implementations of the NGCA and seems to improve significantly the accuracy of estimation.

On the other hand, the proposed algorithm is computationally demanding. While the first-order optimization algorithm we propose allows to treat efficiently the problems which are far beyond the reach of classical SDP-optimization techniques, the numerical difficulty seems to be the main practical limitation of the proposed approach.
%-------------------------------------------------------------------------
\newpage
%\subsubsection*{References}
%-------------------------------------------------------------------------
%{\small

%}
\newpage
\appendix
%-------------------------------------------------------------------------
\section{Appendix }
\label{proofs}
%-------------------------------------------------------------------------

Let $X=X^T\in \R^{L\times L}$ be positive semidefinite with $|X|_1\le 1$, and
let $Y=X^{1/2}$ be the symmetric  positive semidefinite square root of $X$. If we denote
$y_i,\,i=1,..,L$ the columns of $Y$, then $|X|_1\le 1$ implies that
\[
\sum_{1\le i,j\le L}|y_i^Ty_j|\le 1.
\]
We make here one trivial though useful observation: for any matrix
$A\in \R^{d\times L}$, when denoting $B=A^TA$, we have
\be
\|AY\|^2_2&=&\tr(A^TAX)=\tr[BX]=\sum_{j=1}^L \sum_{i=1}^L B_{ji}X_{ij}\le \max_{ij}|B_{ij}|=|A|_2^2.
\ee{frob}
(Recall that for a matrix $A\in \R^{d\times L}$ with columns $a_i$, $i=1,...,L$, $|A|_2$ stands for the maximal column norm:
$
|A|_2=\max_{1\le i\le L}|a_i|_2$).
%In the sequel we also need the ``deterministic analogues'' of %the problem \rf{proj1}
%%\be
%%\min_P\max_c\left[{f}(c,{P})\equiv\left\{\tr\left[\|(I-P){U}c\|\right]\left|\;\begin{array}{c}0\preceq P\preceq I,\;\tr P=m,\; \\
%%\|c\|_1\le 1,\;{G}c=0\end{array}\right.\right\}\right].
%%\ee{opt1}
%%and that of
%the problem \rf{proj3}:
%\be
%\min_P\max_X\left[{f}(X,{P})\equiv\left\{\tr\left[{U}(I-P){U}X\right]\left|\;\begin{array}{c}0\preceq P\preceq I,\;\tr P=m,\; \\
%X\succeq 0,\;|X|_1\le 1,\;\tr[{G}X{G}]=0\end{array}\right.\right\}\right].
%\ee{sdp0}

%Observe that for any $c$ which is feasible for \rf{opt1}, the rank $1$ matrix $X=cc^T$ certainly belongs to the feasible set of \rf{sdp0}.
%So for any feasible $P$ the maximal in $X$ value of the latter problem is not less than the maximum in $c$
%of \rf{opt1}.
%We can rewrite the problems \rf{sdp0} and \rf{proj3} using $Y=X^{1/2}$, so that the objective functions  $g(Y,P)$ of
% \rf{sdp0} and $\hat{g}(Y,P)$ of \rf{proj3} become, respectively,
%\[
%g(Y,P)=\|(I-P)^{1/2}UY\|^2_2,\;\;\;\hat{g}(Y,P)=\|(I-P)^{1/2}\hat{U}Y\|_2^2.
%\]
%Let now $(\hat{X},\hat{P})$ be the saddle point of \rf{proj3}. If we denote $\hat{f}(X,{P})$  the objective of \rf{proj3} then
%\[
%\hat{f}(X,\hat{P})\le[\hat{f}_*\equiv\hat{f}(\hat{X},\hat{P})]\le \hat{f}(\hat{X},P),
%\]
%for any feasible $P$ and $X$. We denote $\hat{Y}=\hat{X}^{1/2}$, $( X^*,P^*)$ (with $Y^*=(X^*)^{1/2}$) the saddle point of \rf{sdp0}.% and $\bar{P}$ an optimal solution to \rf{opt1}.

We can rewrite the problem  \rf{proj3} using $Y=X^{1/2}$, so that the objective
\begin{eqnarray*}\hat{f}(X,P)=\tr[\hat{U}^T(I-P)\hat{U}X]\end{eqnarray*} of \rf{proj3}
becomes
\[
\hat{g}(Y,P)=\|(I-P)^{1/2}\hat{U}Y\|_2^2.
\]
Let now $(\hat{X},\hat{P})$ be a saddle point of \rf{proj3}.
Namely, we have for any feasible $P$ and $X$:
%if $\hat{f}(X,{P})$ is the objective of \rf{proj3} then
\[
\hat{f}(X,\hat{P})\le[\hat{f}_*\equiv\hat{f}(\hat{X},\hat{P})]\le \hat{f}(\hat{X},P),
\]
We denote $\hat{Y}=\hat{X}^{1/2}$.
%We are  now to state some properties

In what follows we suppose that vectors $\gamma_\ell$ and
$\eta_\ell$, $\ell=1,...,L$ satisfy \rf{mdp}. In other words,it holds
$|\hat U-U|_2\le \delta_N$ and $|\hat G-G|_2\le \gamma_N$.
\subsection{Proof of Theorem \ref{the:cond1}.}
\begin{lemma}\label{obj1}
%The saddle-point value  $\hat{f}_*$  of \rf{proj3} satisfies
%\be
%\hat{f}^{1/2}_*\le f^{1/2}_*-\delta_N.
%\ee{lower1}
%where  $f_*$ is the saddle-point value of \rf{opt1}.
%Further, when $m\ge m$,
Let $\hat{P}$ be an optimal solution to  \rf{proj3}. Then
%, then for $\omega \in A$, where $A$ is the set in the definition of $\nu_N$ and $\delta_N$ in \rf{mdp},
\be
\max_c\left\{|(I-\hat{P})^{1/2} Uc|_2\,\left|\;|c|_1\le 1,\;Gc=0\right.\right\}\le \lambda_{\min}^{-1}(\Sigma)(\varrho+\nu_N)+2\delta_N.
\ee{upper1}
\end{lemma}
%\paragraph{Proof of Lemma \ref{obj1}}
%\begin{proof}
{\em Proof.}
%Observe that $X^*=Y^*$ is a feasible solution to the approximated problem \rf{proj3}.
%Indeed,  due to \rf{mdp} we have
%$|\hat{G}-G|_2\le \delta_N$, thus by \rf{frob},
%\[
%\|\hat{G}Y^*\|_2\le \|GY^*\|_2+\|(\hat{G}-G)Y^*\|_2\le |\hat{G}-G|_2\le \delta_N.
%\]
%Now, due to the optimality of $\hat{P}$ and $\hat{Y}$.
%\bse
%\left(\hat{f}_*\right)^{1/2}&=&\|(I-\hat{P})^{1/2}\hat{U}\hat{Y}\|_2\ge \|(I-\hat{P})^{1/2}\hat{U}Y^*\|_2\ge
%\|(I-\hat{P})^{1/2}\hat{U}Y^*\|_2\\
%&\ge &\|(I-\hat{P})^{1/2}UY^*\|_2-|(I-\hat{P})^{1/2}(\hat{U}-U)|_2\ge \|(I-P^*)^{1/2}UY^*\|_2-\delta_N\\
%&\ge& \max_{c} \left\{\|(I-P^*)^{1/2}Uc\|_2\left|\;\|c\|_1\le 1,\;Gc=0\right.\right\}-\delta_N\\
%&\ge& \max_{c} \left\{
%\|(I-\bar{P})^{1/2}Uc\|_2\left|\;\|c\|_1\le 1,\;Gc=0\right.\right\}-\delta_N=(f_*)^{1/2}-\delta_N,
%\ese
%what is \rf{lower1}.
We write:
%  (recall that to $0\preceq I-\hat{P}\preceq I$ and thus $0\preceq (I-\hat{P})^{1/2}\preceq I$):
\bse
\lefteqn{
\max_{c}\left\{|(I-\hat{P})^{1/2}Uc|_2\left|\;|c|_1\le 1,\;Gc=0\right.\right\}
}\\
&\le&
%\max_{x\in \hat{H}}|(I-\hat{P})^{1/2}Ux|\le
\max_{Y} \left\{\|(I-\hat{P})^{1/2}UY\|_2\left|\;\;|Y^2|_1\le 1,\;GY=0\right.\right\}
\\
&\le& \max_{Y} \left\{\|(I-\hat{P})^{1/2}\hat{U}Y\|_2\left|\;\;|Y^2|_1\le 1,\;GY=0\right.\right\}\\
&&+\max_{Y} \left\{\|(I-\hat{P})^{1/2}(\hat{U}-U)Y\|_2\left|\;\;|Y^2|_1\le 1,\;GY=0\right.\right\}\\
&\le&\max_{Y}
\left\{\|(I-\hat{P})^{1/2}\hat{U}Y\|_2\left|\;\;|Y^2|_1\le 1,\;\|\hat{G}Y\|_2\le \varrho\right.\right\}\\
\mbox{(by \rf{frob})}&&+|(I-\hat{P})^{1/2}(\hat{U}-U)|_2\\
\mbox{(due to  $0\preceq I-\hat{P}\preceq I$)}&= &\|(I-\hat{P})^{1/2}\hat{U}\hat{Y}\|_2+\delta_N\le
 \|(I-\Pi^*)^{1/2}\hat{U}\hat{Y}\|_2+\delta_N
\\\mbox{(again by \rf{frob})}&\le&
%\|(I-P^*)^{1/2}U\hat{Y}\|_2+2\delta_N\le
\|(I-\Pi^*)^{1/2}U\hat{Y}\|_2+2\delta_N.
%\\
%&\le& |\Sigma^{-1}|\,\|\hat{G}-G\|+3\delta_N\le C(\Sigma)\delta_N+3\delta_N.
\ese
On the other hand, as $\|\hat{G}\hat{Y}\|_2\le \nu_N$, we get
\[
\|G\hat{Y}\|_2\le \|\hat{G}\hat{Y}\|_2+\|(\hat{G}-G)\hat{Y}\|_2\le \varrho+|\hat{G}-G|_2\le \varrho+\nu_N,
\]
and by \rf{bounds},
\[
\|(I-\Pi^*)U\hat{Y}\|_2\le \lambda_{\min}^{-1}(\Sigma)(\varrho+\nu_N).
\]
 This implies \rf{upper1}.
\epr
%  \begin{proof}
We now come back to the proof of the theorem. Let $\hat\lambda_j$ and $\hat\theta_j$, $j=1,\ldots,d$ be
respectively the eigenvalues and the eigenvectors of $\widehat{P}$. Assume
that $\hat\lambda_1\ge \hat\lambda_2\ge \ldots\ge \hat\lambda_d$. Then
$\widehat{P}=\sum_{j=1}^d\hat\lambda_j \hat\theta_j\hat\theta_j^{T}$ and
$\widehat\Pi=\sum_{j=1}^{m}\hat\theta_j\hat\theta_j^{T}$. Let $\beta=Uc$ for $c$ such that $|c|_1\le 1$ and $Gc=0$. We have
 \bse
 \beta^T(I-\hat P)\beta&=&\sum_{j=1}^m (1-\hat\lambda_j)(\hat\theta_j^{T}\beta)^2+
 \sum_{j>m}(1-\hat\lambda_j)(\hat\theta_j^{T}\beta)^2\\
 &\ge&   \sum_{j>m}(1-\hat\lambda_j)(\hat\theta_j^{T}\beta)^2\ge (1-\hat\lambda_{m+1})(\hat\theta_j^{T}\beta)^2\\
 &=&(1-\hat\lambda_{m+1})\beta^T(I-\hat\Pi)\beta=(1-\hat\lambda_{m+1})|(I-\hat\Pi)\beta|^2_2.
 \ese
Since, for obvious reasons, $\hat\lambda_{m+1}\le {m\over m+1}$, it applies (i) due to \rf{upper1}.

Let us show (ii). We have due to \rf{ident} and \rf{upper1}:
 \be
\lefteqn{ \tr\left[(I-\hat{P})\Pi^*\right]=
 \tr\left[(I-\hat{P})^{1/2}\Pi^*(I-\hat{P})^{1/2}\right]}\nn
 &\le& \sum_{k=1}^{{m}}\mu_k\tr\left[(I-\hat{P})^{1/2}Uc_kc_k^TU^T(I-\hat{P})^{1/2}\right]
 =\sum_{k=1}^{{m}}\mu_k|(I-\hat{P})^{1/2}Uc_k|_2^2 \nn
 &\le&
  \sum_{k=1}^{{m}}\mu_k \max_{c}\left\{|(I-\hat{P})^{1/2}Uc|_2^2\left|\;|c|_1\le 1,\;Gc=0\right.\right\}\nn
&  =&\mu^*(\lambda_{\min}^{-1}(\Sigma)(\varrho+\nu_N)+2\delta_N)^2,
\ee{ima1}
which is \rf{id1}.\\[2ex]

 %Therefore, \bse
%\tr[\widehat{P}\Pi^*]&\le& \sum_{j\le m} \hat\lambda_j \hat\theta_j^{T}
%\Pi^*\hat\theta_j+
%\hat\lambda_{m+1}\sum_{j> m}\hat\theta_j^{T} \Pi^*\hat\theta_j\\
%&=&\sum_{j\le m} \hat\lambda_j\hat\theta_j^{T} \Pi^*\hat\theta_j+
%\hat\lambda_{m+1}\tr[(I-\hat\Pi)\Pi^*]\\
%&=&\sum_{j\le m} (\hat\lambda_j-\hat\lambda_{m+1})
%\hat\theta_j^{T} \Pi^*\hat\theta_j+m\hat\lambda_{m+1}\\
%\mbox{(since $\hat\theta_j^{T} \Pi^*\hat\theta_j=|\Pi^*\hat\theta_j|^2\le 1$)}&\le &\sum_{j\le m} \hat\lambda_j.
%\ese
% Now
Note that $\tr[\widehat{P}\Pi^*]\le \sum_{j\le m} \hat\lambda_j$ (cf, e.g., Corollary 4.3.18 of \cite{HornJ}), thus
by \rf{ima1},
\[
\hat\lambda_{m+1}\le m-\sum_{j\le m}\hat\lambda_j \le
\tr[(I-\widehat{P})\Pi^*]\le\mu^*(\lambda_{\min}^{-1}(\Sigma)(\varrho+\nu_N)+2\delta_N)^2.
\]
On the other hand,
 \bse
 \tr[(I-\hat P)\Pi^*]&=&\sum_{j=1}^m (1-\hat\lambda_j)\hat\theta_j^{T} \Pi^*\hat\theta_j+
 \sum_{j>m}(1-\hat\lambda_j)\hat\theta_j^{T} \Pi^*\hat\theta_j\\
 &\ge&   \sum_{j>m}(1-\hat\lambda_j)\hat\theta_j^{T} \Pi^*\hat\theta_j\ge (1-\hat\lambda_{m+1})\hat\theta_j^{T} \Pi^*\hat\theta_j\\
 &=&(1-\hat\lambda_{m+1})\tr[(I-\hat\Pi)\Pi^*],
 \ese
 and we conclude that
\[
\tr[(I-\widehat{\Pi})\Pi^*]\le {\tr[(1-\hat P)\Pi^*]\over 1-\hat\lambda_{m+1}}\le {\mu^*(\lambda_{\min}^{-1}(\Sigma)(\varrho+\nu_N)+2\delta_N)^2\over 1-\mu^*(\lambda_{\min}^{-1}(\Sigma)(\varrho+\nu_N)+2\delta_N)^2}.
\]
Now, using the relation $\tr \hat{\Pi}=\tr \Pi^*=m$, we come to
\[
\|\hat{\Pi}-\Pi^*\|_2^2=\tr[\hat{\Pi}^2-2\hat{\Pi}\Pi^*+(\Pi^*)^2]
= 2m-2\tr[\hat{\Pi}\Pi^*]=
2\tr[(I-\hat{\Pi})\Pi^*],
\]
and we arrive at \rf{id2}.
\subsection{Proof of Theorem \ref{the:cond2}.}
Let now  $\hat{P}, \;\hat X $ and $\hat{t}=\tr \hat P$ be a triplet of optimal solution to  \rf{proj4}.
\begin{lemma}\label{obj2}
Let $\hat{P}$ be an optimal solution to  \rf{proj4}.
\item[(i)] In the premises of the theorem $\Pi^*$ is a feasible solution of \rf{proj4} and
$\tr \hat P\le \tr \Pi^*=m$.
\item[(ii)] We have
\be
\max_c\left\{|(I-\hat{P})^{1/2} Uc|_2\,\left|\;|c|_1\le 1,\;Gc=0\right.\right\}\le \rho+\delta_N.
\ee{upper12}
\end{lemma}
%\paragraph{Proof of Lemma \ref{obj1}}
%\begin{proof}
{\em Proof.}
We act as in the proof of Lemma \ref{obj1}: to verify $(i)$ we observe that
\bse
\lefteqn{\max_{X}\left\{\tr[\hat U^T(I-\Pi^*)\hat U X]\left|\;X\succeq 0,\;|X|_1\le 1,\;\tr[\hat GX\hat G ^T]\le\varrho^2\right.\right\}}\\
&=&\max_{Y}\left\{
\|(I-\Pi^*)\hat U Y\|^2_2\left|\;\;|Y^2|_1\le 1,\;\|\hat GY\|_2\le\varrho\right.\right\}\\
&\le& \max_{Y}\left\{
(\|(I-\Pi^*)U Y\|_2+\delta_N)^2\left|\;\;|Y^2|_1\le 1,\;\|\hat GY\|_2\le\varrho\right.\right\}\\
&\le &
\left(\lambda^{-1}_{min}(\Sigma)(\varrho+\nu_N)+\delta_N\right)^2.
\ese
Thus, if $\rho\ge \lambda^{-1}_{min}(\Sigma)(\varrho+\nu_N)+\delta_N$, $\Pi^*$ is a feasible solution of \rf{proj4} and, as a result, $\tr \hat P\le \tr \Pi^*$.
To show $(ii)$ it suffices to note that
\bse
\lefteqn{
\max_{c}\left\{|(I-\hat{P})^{1/2}Uc|_2\left|\;|c|_1\le 1,\;Gc=0\right.\right\}
}\\
&\le&
\max_{Y} \left\{\|(I-\hat{P})^{1/2}UY\|_2\left|\;\;|Y^2|_1\le 1,\;GY=0\right.\right\}
\\
&\le& \max_{Y} \left\{\|(I-\hat{P})^{1/2}\hat{U}Y\|_2\left|\;\;|Y^2|_1\le 1,\;GY=0\right.\right\}+|(I-\hat{P})^{1/2}(\hat{U}-U)|_2\\
&\le& \max_{Y} \left\{\|(I-\hat{P})^{1/2}\hat{U}Y\|_2\left|\;\;|Y^2|_1\le
1,\;\|\hat GY\|_2\le \varrho\right.\right\}+\delta_N \le \rho+\delta_N
\ese because of the feasibility of $\hat P$. \epr
Now using the bound $\hat m\le m$ we complete the proof
following exactly the lines of the proof of Theorem \ref{the:cond1}.

\subsection{Proof of Proposition \ref{arik100}}
\label{sec:arik}
Observe that
\be
\underline{F}( \hat Z, \hat Y) &=& \min_{(P,W)\in \cP\times \cW}
\left\{
\tr[B^T (I - P)B Q^T \hat Z Q] + \lambda\, \tr[W(Q^T \hat Z Q - \hat Y )]\right\}\nn
&=&\min_{P\in \cP}
\left\{
\tr[B^T (I - P)B Q^T \hat Z Q] - \lambda\, \|Q^T \hat Z Q - \hat Y \|_2\right\}\nn
&\le& \min_{P\in \cP}\left\{
\tr[B^T (I - P)B Q^T \hat ZQ] \right\}
= \underline{G}(\hat Z);
\ee{5}
and
\bse
\bar{F}(\hat P,\hat W) &=& \max_{(Z,Y )\in \cZ\times \cY}
\left\{\tr[B^T(I -\hat P)B Q^TZQ] + \lambda\, \tr[\hat W (Q^TZQ - Y)]\right\}\\
&\ge &\max_{Z\in \cZ,\;|Q^TZQ|_1\le 1,\,Y=Q^TZQ}
\left\{\tr[B^T(I -\hat P)B Q^TZQ] + \lambda\, \tr[\hat W (Q^TZQ - Y)]\right\}\\
&=& \bar{G}(\hat P):
\ese
Assume first that $|Q^T \hat ZQ|_1\le 1$. In this case $\tilde Z = \hat Z$  and
\[ \e\ge \bar{F}(\hat P,\hat W) - \underline{F}(\hat Z, \hat Y) = \bar{G}(\hat P)-
\underline{G}(\hat Z)=\bar{G}(\hat P)-
\underline{G}(\tilde{Z})
\]
(the second ¸ is given by \rf{5}), as claimed.
Now assume that $s=Q^T \hat Z Q|_1 > 1$. We have already established the first equality
of the following chain:
\bse
\lefteqn{\underline{F}( \hat Z, \hat Y ) = \min_{P\in \cP}
\left\{\tr(B^T (I- P)BQ^T \hat ZQ)-\lambda\, \|Q^T \hat Z Q -\hat Y \|_2\right\}
}\\
&\le&
\min_{P\in \cP}
\left\{\tr[B^T (I- P)BQ^T \hat ZQ]-{\lambda\over L} |Q^T \hat Z Q -\hat Y |_1\right\}\\
&\le &\min_{P\in \cP}
\left\{\tr[B^T (I- P)BQ^T \hat ZQ]-{\lambda\over L} (s -1)\right\}\\
&=& \min_{P\in \cP}
\left\{s\tr[B^T (I- P)BQ^T \tilde{Z}Q]-{\lambda\over L} (s -1)\right\}\\
&\le &\min_{P\in \cP}
\left\{\tr[B^T (I- P)BQ^T \tilde{Z}Q]+
\underbrace{(s-1)|B^T(I-P)B|_\infty|Q^T\tilde{Z}Q|_1}_{\le (s-1)|B^TB|_\infty=(s-1)|B|_2^2}
-{\lambda\over L} (s -1)\right\}\\
&\le& \min_{P\in \cP}
\left\{\tr[B^T (I- P)BQ^T \tilde{Z}Q]\right\}=\underline{G}( \tilde{Z}),
\ese
where the concluding $\le $ is readily given by the definition of $\lambda$.\footnote{We denote $|A|_\infty=\max_{ij}|A_{ij}|$.} Further, we have already
seen that
\[
\bar{F} (\hat P,\hat W)\ge \bar{G}(\hat P).
\]
Consequently,
\[
\e\ge \bar{F}(\hat P,\hat W)- \underline{F}(\hat Z, \hat Y )\ge  \bar{G}(\hat P)-\underline{G}(\tilde{Z}),
\]
as claimed.\epr
\subsection{Computing the prox-transform}
\label{sec:prox}
Recall that because of the additivity of the distance-generating function $d$ the computation of the prox-transform on the set $\cP\times \cW \times \cZ\times \cY$ can be decomposed into independent computations on the four domains of \rf{2}.

\paragraph{Prox-transform on $\cP$.}
The proxy-function of $\cP$ is the matrix entropy:
\begin{eqnarray*}d(P_0,P)=\beta_P\tr\left[{P\over m}\left(\ln \left({P\over
m}\right)-\ln \left({P_0\over m}\right)\right)\right]\;\;\mbox{for}\;\; P,P_0\in
\cP,\;\;\beta_P>0.\end{eqnarray*} To compute the corresponding component of $T$ we need
to find, given $S\in {\bold S}^d$,
\begin{eqnarray}
\label{eq:TP}
T_\beta(P_0,S)&=&\arg\max_{P\in \cal P} \left\{\tr[S(P-P_0)]-\beta_P \,\tr\left[{P\over m}\left(\ln \left({P\over m}\right)-\ln \left({P_0\over m}\right)\right)\right]\right\}\\
&=& \arg\max_{P\in \cal P} \left\{\tr\left[\left(S+{\beta_P\over m} \ln \left({P_0\over m}\right)\right)P\right]-\beta_P \,\tr\left[{P\over m}\ln \left({P\over m}\right)\right]\right\}.\nonumber
\end{eqnarray}
By the symmetry considerations we conclude that the optimal solution of this problem is
diagonal in the basis of eigenvectors of $ S+{\beta_P\over m} \ln ({P_0\over m})$. Thus
the solution of  \rf{TP} can be obtained as follows: compute  the eigenvalue
decomposition  \begin{eqnarray*}S+{\beta_P\over m} \ln ({P_0\over m})=\Gamma
\Lambda\Gamma^T\end{eqnarray*} and let $\lambda$ be the diagonal of $\Lambda$. Then
solve the ``vector" problem \be p^*=\arg\max_{0\le p\le 1,\;\sum p\le m}
\lambda^Tp-{\beta_P\over m} \sum_{i=1}^d p_i\ln(p_i/m). \ee{ent1} and compose
\[
T_\beta(P,S)=\Gamma {\rm diag}(y^*)\Gamma^T.
\]
\par
Now, the solution of \rf{ent1} can  be obtained by simple bisection: indeed, using Lagrange duality we conclude that the components of $y^*$ satisfies
\[
p^*_i=\exp\left({\lambda_i\over \beta}-\nu\right)\wedge 1,\;\;i=1,...,d,
\]
 and the Lagrange multiplier $\nu$ is to be set to obtain $\sum p^*_i=m$, what can be done by bisection in $\nu$. When the solution is obtained, the optimal value of \rf{TP} can be easily computed.

\paragraph{Prox-transform on $\cW$.} The distance-generating function of $\cW$ is  $\beta_W\tr[W^2]/2=\|W\|_2^2/2$ so that we have to solve for $S\in {\bold S}^L$
\begin{eqnarray}
\label{eq:TW}
T_\beta(W_0,S)&=&\arg\max_{\|W\|_2\le 1} \left\{\tr[S(W-W_0)]-{\beta_W\over 2} {\|W-W_0\|_2^2\over 2}\right\}.
\end{eqnarray}
The optimal solution to \rf{TW} can be easily computed
\[
T_\beta(W_0,S)=\left\{\begin{array}{lcl}W_0+S/\beta_W&\mbox{if}&\|W_0+S/\beta_W\|_2\le 1,\\
(W_0+S/\beta_W)/\|W_0+S/\beta_W\|_2&\mbox{if}&\|W_0+S/\beta_W\|_2> 1.
\end{array}\right.
\]

\paragraph{Prox-transform on $\cZ$.} The prox-function of $\cZ$ is  the matrix entropy and we have to solve for $S\in {\bold S}^{L-d}$
\bse
T_\beta(Z,S)&=&\arg\max_{Z\in \cal Z} \tr[S(Z-Z_0)]-\beta_Z \tr[Z(\ln (Z)-\ln (Z_0))]\\
&=& \arg\max_{Z\in \cal Z} \tr[(S+\beta_Z \ln (Z_0))Z]-\beta_Z \tr[Z\ln Z].\nonumber
\ese
Once again,  in the  basis of eigenvectors of $ S+\beta_Z \ln (Z_0)$ the problem reduces to
\bse
z^*=\arg\max_{z\ge 0,\;\sum z\le 1} \lambda^Tz-{\beta_Z} \sum_{i=1}^d z_i\ln(z_i),
\ese%{ent2}
where $\lambda$ is the diagonal of $\Lambda$ with  $S+\beta_Z \ln Z_0=\Gamma \Lambda\Gamma^T$.
In this case
\[
z^*_i=\frac{\exp({\lambda_i\over \beta})}{\sum _{j=1}^{\bar{L}}
\exp({\lambda_j\over \beta})},\;\;i=1,...,L-d.
\]

\paragraph{Prox-transform on $\cY$.} The distance generating function for the domain $\cY$ is defined as follows:
\bse
d(Y)  &=&  \min \left\{\;  \sum\limits_{i,j=1}^L \left(
u_{ij}\ln[u_{ij}]  +  v_{ij}\ln[v_{ij}] \,\right): \;
\sum\limits_{i=1}^n ( u_{ij}+ v_{ij} ) = 1,\right.\\
 Y_{ij} &=& \left. u_{ij} - v_{ij},\; u_{ij} \geq 0,\;
v_{ij} \geq 0,  \; 1\le i,j \le L \; \right\}.% +\ln [2L^2],
%\\
%\\
%\psi(t ) &=& \left\{ \begin{array}{rl} t \ln t, & t > 0, \\ 0, &
%t = 0 . \end{array} \right.
\ese
In other words, the element $Y\in \cY$ is decomposed according to $Y=u-v$, where $(u,v)$ is an element of the $2L^2$-dimensional simplex $\Delta=\left\{x\in \R^{2L^2},\,x\ge 0,\, \sum_i x_i=1\right\}$. To find the $Y$-component of the prox-transform amounts to find for $S\in {\bold S}^L$
\be
T_{\beta_Y}(Y_0,S)&=&T_{\beta_Y}(u^0,v^0,S)\nn&=&\arg\max_{u,v\in \Delta}
\tr[S(u-v)]-\beta_Y \sum_{ij} \left[u_{ij}\ln ({u_{ij}\over u^0_{ij}})
+v_{ij}\ln ({v_{ij}\over v^0_{ij}})\right].
\ee{TY}
One can easily obtain an explicit solution to \rf{TY}: let
\begin{eqnarray*}
a_{ij}=u^0_{ij}\exp\left({S_{ij}\over \beta_Y}\right),\;\;\;
b_{ij}=v^0_{ij}\exp\left(-{S_{ij}\over \beta_Y}\right).
\end{eqnarray*}
Then $T_{\beta_Y}(Y_0,S)=u^*-v^*$, where
\[
u^*_{ij}={a_{ij}\over \sum_{ij} (a_{ij}+b_{ij})},\;\;v^*_{ij}={b_{ij}\over \sum_{ij} (a_{ij}+b_{ij})}.
\]
\end{document}